\documentclass[10pt,a4paper]{article}
\usepackage{amsmath, amssymb}
\usepackage{latexsym, color}
\usepackage{epsfig}

\begin{document}

\newtheorem{thm}{Theorem}[section]
\newtheorem{cor}[thm]{Corollary}
\newtheorem{lem}[thm]{Lemma}
\newtheorem{prop}[thm]{Proposition}
\newtheorem{defn}[thm]{Definition}
\newtheorem{rem}[thm]{Remark}
\newtheorem{Ex}[thm]{EXAMPLE}
\def\nm{\noalign{\medskip}}

\bibliographystyle{plain}

%\numberwithin{equation}{section}

\newcommand{\qed}{\hfill \ensuremath{\square}}
\newcommand{\ds}{\displaystyle}
\newcommand{\pf}{\medskip \noindent {\sl Proof}. ~ }
\newcommand{\p}{\partial}
\renewcommand{\a}{\alpha}
\newcommand{\z}{\zeta}
\newcommand{\pd}[2]{\frac {\p #1}{\p #2}}
\newcommand{\norm}[1]{\| #1 \|}
\newcommand{\dbar}{\overline \p}
\newcommand{\eqnref}[1]{(\ref {#1})}
\newcommand{\na}{\nabla}
\newcommand{\Om}{\Omega}
\newcommand{\ep}{\epsilon}
\newcommand{\tmu}{\widetilde \mu}
\newcommand{\vep}{\varepsilon}
\newcommand{\tlambda}{\widetilde \lambda}
\newcommand{\tnu}{\widetilde \nu}
\newcommand{\vp}{\varphi}
\newcommand{\RR}{\mathbb{R}}
\newcommand{\CC}{\mathbb{C}}
\newcommand{\NN}{\mathbb{N}}
\renewcommand{\div}{\mbox{div}~}
\newcommand{\bu}{{\bf u}}
\newcommand{\la}{\langle}
\newcommand{\ra}{\rangle}
\newcommand{\Scal}{\mathcal{S}}
\newcommand{\Lcal}{\mathcal{L}}
\newcommand{\Kcal}{\mathcal{K}}
\newcommand{\Dcal}{\mathcal{D}}
\newcommand{\tScal}{\widetilde{\mathcal{S}}}
\newcommand{\tKcal}{\widetilde{\mathcal{K}}}
\newcommand{\Pcal}{\mathcal{P}}
\newcommand{\Qcal}{\mathcal{Q}}
\newcommand{\id}{\mbox{Id}}
%%%%%%%%%%
\newcommand{\be}{\begin{equation}}
\newcommand{\ee}{\end{equation}}

\title{Strong Influence of a Small Fiber on Shear Stress in Fiber-Reinforced Composites}

\author{Mikyoung Lim\thanks{Department of Mathematics, Colorado State
University, Fort Collins, CO 80523, USA (lim@math.colostate.edu)}
\and KiHyun Yun\thanks{\footnotesize Department of Mathematics,
Michigan State University, East Lansing, MI 48824,
USA (kyun@math.msu.edu)}}
 \maketitle

\begin{abstract}
In stiff fiber-reinforced material, the high shear stress
concentration occurs in the narrow region between fibers. With the
addition of a small geometric change in cross-section, such as a
thin fiber or a overhanging part of fiber, the concentration is
significantly increased. This paper presents mathematical analysis
to explain the rapidly increased growth of the stress by a small
particle in cross-section. To do so, we consider two crucial cases
where a thin fiber exists between a pair of fibers, and where one of
two fibers has a protruding small lump in cross-section. For each
case, the optimal lower and upper bounds on the stress associated
with the geometrical factors of fibers is established to explain the
strongly increased growth of the stress by a small particle.
\end{abstract}

MSC-class:  35J25, 73C40

\section{Introduction}
\begin{figure}[h!]
\begin{center}
\epsfig{figure=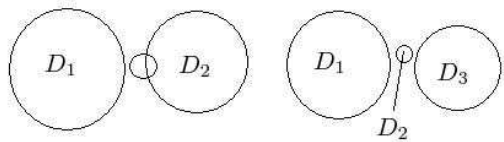,width=7.0cm}
\end{center}
\caption{Case (A) and case (B)}\label{caseAB}\end{figure}

\begin{figure}[h!]
\begin{center}\includegraphics[width=4 cm]{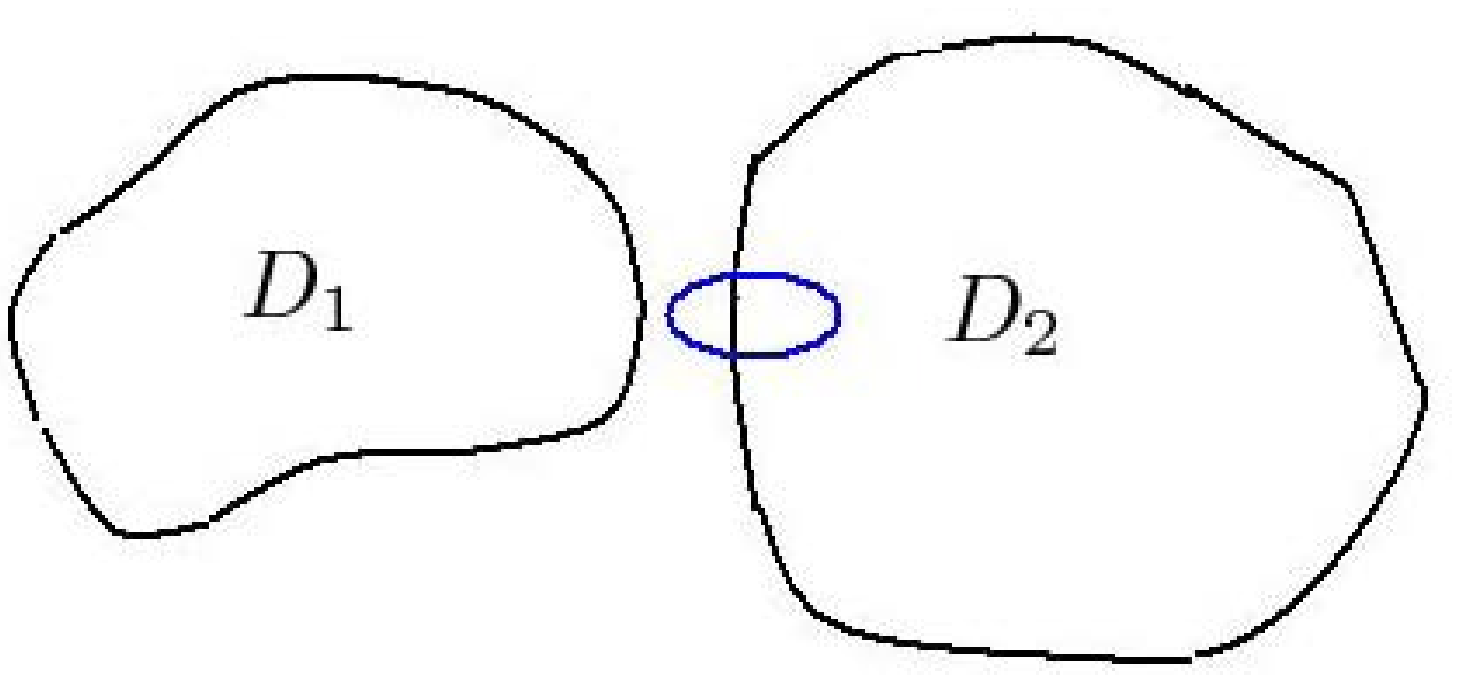}
%\hbox{\setlength{\epsfxsize}{3.2 in} \epsfbox{96-emmer.eps}}
%\caption{a soup bubble in the tetrahedron}\label{fig.1}
~~~\includegraphics[width=4 cm]{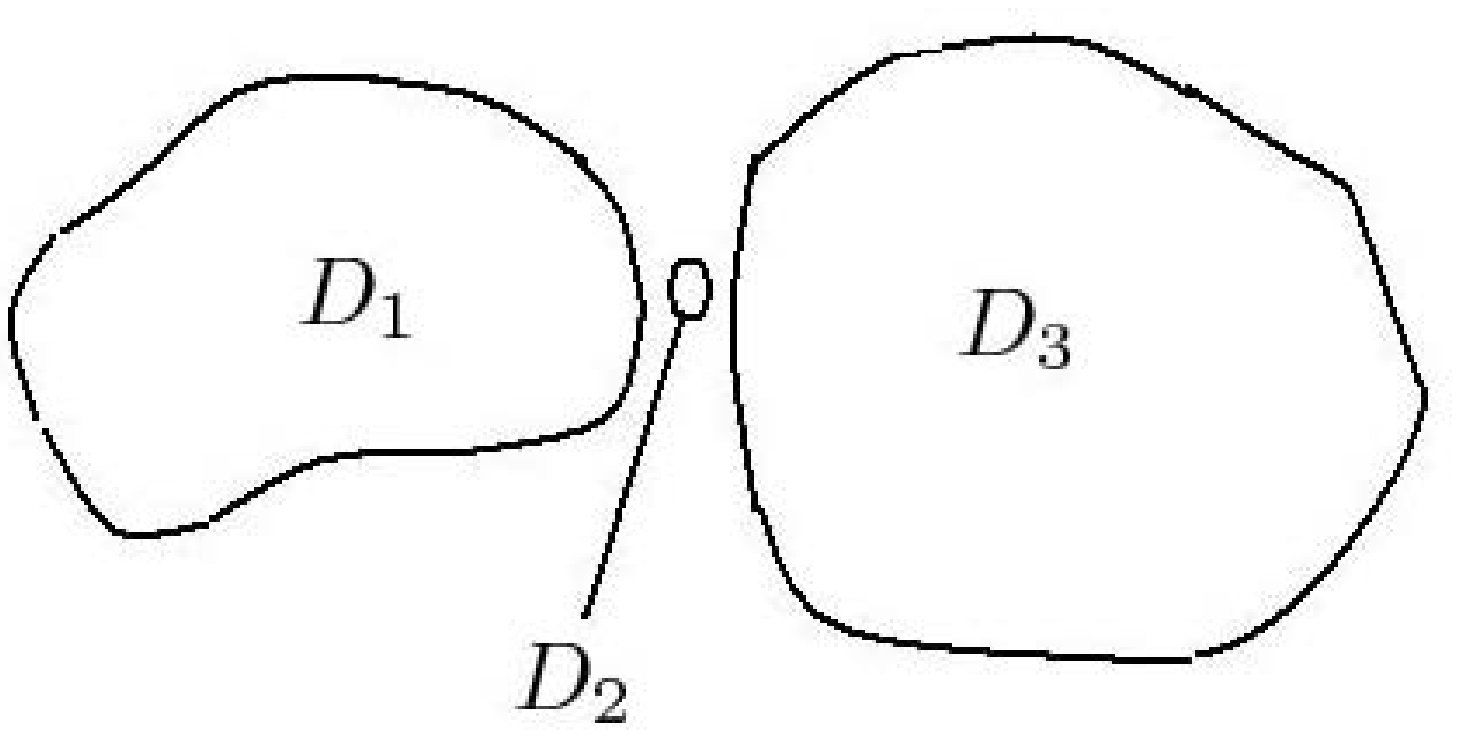}
%\hbox{\setlength{\epsfxsize}{3.2 in} \epsfbox{96-emmer.eps}}
%\caption{caption Caption Caption Caption.} \label{fig.1}
\caption{Case (C) and case (D)}\label{caseCD}
\end{center}
\end{figure}

In this paper, we concern ourselves with the high stress
concentration occurring in the stiff fiber-reinforced composites
when fibers are located closely. The primary investigation focuses
on the case when a smaller fiber is located in-between area of two
fibers, see Figure \ref{caseAB} and Figure \ref{caseCD}. This paper
reveals that, with the addition of a smaller fiber, the growth of
stress is significantly increased: if the diameter $d$ of the fiber
in the middle is sufficiently small and the distance between
adjoining fibers is $\epsilon$, then the stress blows up at the rate
of $\frac 1 {\sqrt{d \epsilon}}$ in the narrow region, even though
the blow-up rate has been known as $\frac 1 {\sqrt {\epsilon}}$ as
in the case of a pair of fibers. This means that the defect of fiber
as a protrusion causes much lower strengths in composites than had
been thought. To derive it, we estimate the optimal lower and upper
bounds of the stress concentration in terms of the diameters of
fibers and the distances between them. These bounds explain the
dramatic change of the growth of stress when the diameter of the
fiber placed in the middle is relatively smaller than other two
fibers.

\par In the anti-plane shear model, the stress tensor represents the
electric field in the two dimensional space, where the out-of-plane
elastic displacement satisfies a conductivity equation, and the
cross-section of stiff fibers corresponds to the embedded
conductors. In this respect, we consider the gradient of the
solution to a conductivity problem to estimate the stress. Adjacent
stiff fiber-reinforcement induces the high stress concentration in
the narrow region between the fibers. This implies the blow-up of
the gradient of the solution between adjoining conductors, see
\cite{ADKL, AKL, AKLLL, BC, BLY, K,Y,Y2}.

\par  Meanwhile, the extreme conductivities are indispensable to the
blow-up phenomena: when the inclusion's conductivity is away from
zero and infinity, the boundedness of the stress function has been
derived by Li and Vogelius \cite{LV}, see also \cite{BV}, and it was
generalized to elliptic systems by Li and Nirenberg \cite{LN}. In
\cite{AKL,AKLLL}, for conductivities including both bounded and
extreme cases, Ammari et al.
 have established the optimal bounds of the gradient of solutions to the conductivity equation,
 when conductors are of circular shape in two dimensions, and the optimal bounds provides $\epsilon^{-1/2}$
  blow-up rate, where $\epsilon$ is the distance between two conductors. Yun \cite{Y, Y2} has extended
  this blow-up result for the case of two adjacent perfect conductors of a sufficiently general shape
   in two dimensions. In Bao, Li and Yin's paper \cite{ BLY},
   it has been also investigated as the blow-up phenomena in higher dimensional spaces, also see \cite{ADKL, LY}.
They has also done a natural follow-up in \cite{Bao_thesis, BLY2}
that the blow-up rate known only for a pair of fibers is still valid
for the multiple inclusions in any dimensions.
\par In contrast, our paper witnesses an unexpected fact on multiple inclusions that the
growth of stress can be significantly increased by a little
geometric change of an inclusion, even though the blow-up rate is
still $\epsilon^{-1/2}$.

For $l=1,\dots,L$, let $D_l$ be conducting inclusions in $\RR^2$,
 that is cross-sections of stiff fibers. Then, under the action of the applied field $H$, the electric potential $u$ satisfies the following conductivity equation:
\begin{equation} \label{eq:conductivity}
\quad \left\{
\begin{array}{ll}
\ds\Delta u  = 0,\quad&\mbox{in }{\mathbb{R}^2 \backslash \overline{\cup_{l=1}^L D_l}},\\
\ds u(\textbf{x})- H(\textbf{x}) = O(|\textbf{x}|^{-1}),\quad&\mbox{as } |\textbf{x}| \rightarrow \infty,\\
\ds u |_{\partial D_l} = C_l~ \mbox{(constant)},\quad&\mbox{for }l=1,\dots,L,\\
 \int_{\partial D_l} {\partial_{  \nu} u }~dS = 0,\quad
&\mbox{for } l =1,\dots,L,
\end{array}
\right.
\end{equation}
where $H$ is an entire harmonic function $H$ in $\mathbb{R}^2$ and
$\textbf{x} =(x_1,x_2)$. In this paper, we only consider the case of
$L=2,3$,

As we mentioned previously, for the two fibers with the circular
cross-sectional shape, Ammari et al. \cite{AKL, AKLLL} obtained the
optimal blow-up rate $\epsilon^{-1/2}$ for $\nabla u$, and this
result is extended by Yun \cite{Y, Y2} to general shaped fibers.
Building on the prior results, we extend these into the the
following interesting direction: we first consider the circular
inclusions in Case (A) and Case (B), and second extend the result
into general shaped ones in Case (C) and Case (D). In Case (A) and
Case (B), we add a small circular inclusion between two others so
that three disk centers are lined up in one straight line. The
additional disk can be embedded disjointly from other disks, or it
partially overlap one of two disks, and we formulate these two cases
as follows.
\begin{itemize}
 \item[(A)] One disk and a pair of partially overlapping disks, Figure \ref{caseAB}: there is a portion of disk protruding from one of circular inclusions, i.e.,
$L=2$, and $D_1$ and $D_2$ are $\epsilon$-distanced
domains defined as
\begin{equation}\label{def:proballs}D_1=B_{r_1}(\mathbf{c}_1)\mbox{ and }D_2=B_{r_2}(\mathbf{c}_2)\cup B_{r_3}(\mathbf{c}_3),
\end{equation}
where $B_{r_l}(\mathbf{c}_l)$ is the disk with the radius $r_i$ and centered at $\mathbf{c}_l$,
and \begin{equation}\label{eqn:lumpcenters}\mathbf{c}_1=(-r_1-\frac{\epsilon}{2},0),\ \mathbf{c}_2=(r_2+\frac{\epsilon}{2},0),\mbox{ and } \mathbf{c}_3=(r_3+a+\frac{\epsilon}{2},0).\end{equation}
Here, $B_{r_2}(\mathbf{c}_2)$ is a small disk protruding from $B_{r_3}(\mathbf{c}_3)$, and we assume
$$
B_{r_2}(\mathbf{c}_2)\cap B_{r_3}(\mathbf{c}_3)\neq\emptyset,\mbox{
i.e., }0<a<2r_2,$$  $$\mbox{dist} (B_{r_1}(\mathbf{c}_1),
B_{r_3}(\mathbf{c}_3)) \simeq r_2$$ and $$0<\epsilon\ll r_2\ll
r_1\simeq\ r_3.$$

\item[(B)] Three disjoint disks, Figure \ref{caseAB}: a small disk is disjointly embedded into the in-between area of two disks, i.e.,
$L=3$, and
\begin{equation}\label{def:Bthree}D_l = B_{r_l}(\mathbf{c}_l),\
l=1,2,3,\end{equation}  where $\mathbf{c}_1=(-r_1 - \frac
{\epsilon_1} 2, 0)$, $\mathbf{c}_2=(r_2+ \frac {\epsilon_1} 2,0)$
and $\mathbf{c}_3=(r_3+r_2+\frac {\epsilon_1} 2 + \epsilon_2,0)$.
Hence, the distance between $D_1$ and $D_2$ is $\epsilon_1$, and the
distance between $D_2$ and $D_3$ is $\epsilon_2$. Here, $D_2$ is
regarded as the cross-section of the thin fiber between a pair of fibers with the cross-section $D_1$ and $D_3$. Thus, we assume
that
$$0<\epsilon_{i}\ll r_2\ll r_1\simeq r_3~~\mbox{for}~i =1,2.$$
\end{itemize}

In both cases, the blow-up rate is remarkably increased due to the
existence of $B_{r_2}(\mathbf{c_2})$ as follows:
\begin{thm}[Case A: Protruding small disk]\label{thm:A}
Let $D_1$ and $D_2$ be defined as \eqref{def:proballs}. Then there
is a positive constant $C$ independent of $\epsilon$, $r_1$, $r_2$
and $r_3$ such that
$$u\Bigr|_{\p D_2}-u\Bigr|_{\p D_1}\geq C\frac{r_1r_3}{r_1+r_3}\frac{1}{\sqrt{r_2}}\sqrt{\epsilon},$$
where $u$ is the solution to \eqref{eq:conductivity} with $H(x_1,x_2)=x_1$.
As a result, by the Mean Value Theorem, there is a point $\mathbf{x}_0$ in the narrow region between $D_1$ and $D_2$ such that
$$|\nabla u (\mathbf{x}_0)|\geq C\frac{r_1r_3}{r_1+r_3}\frac{1}{\sqrt{r_2}}\frac 1 {\sqrt{\epsilon}}.$$
\par For any entire harmonic function $H$, let $u$ be the solution
to \eqref{eq:conductivity} with $H$. Then, there is a positive
constant $C$ independent of $\epsilon$, $r_1$, $r_2$ and $r_3$ such
that
$$|\nabla u (\mathbf{x})|\leq C\frac{r_1r_3}{r_1+r_3}\frac{1}{\sqrt{r_2}}\frac 1 {\sqrt{\epsilon}}
$$ in the narrow region between $D_1~\mbox{and}~D_2$.

\end{thm}

\begin{thm}[Case B:  Disjointly embedded small disk] \label{thm:B}
Let $D_i$, $i=1,2,3$, be balls defined as \eqref{def:Bthree}. Then
there is a positive constant $C$ independent of $\epsilon_1$,
$\epsilon_2$, $r_1$, $r_2$ and $r_3$ such that
$$u\Bigr|_{\p D_2}-u\Bigr|_{\p D_1}\geq C\frac{r_1r_3}{r_1+r_3}\frac{1}{\sqrt{r_2}}\sqrt{\epsilon_1},$$
and
$$u\Bigr|_{\p D_3}-u\Bigr|_{\p D_2}\geq C\frac{r_1r_3}{r_1+r_3}\frac{1}{\sqrt{r_2}}\sqrt{\epsilon_2},$$
where $u$ is the solution to \eqref{eq:conductivity} with $H(x_1,x_2)=x_1$. As a result, by the Mean Value Theorem, there exists points; $\mathbf{x}_1$ in the narrow region between $D_1$ and $D_2$; $\mathbf{x}_2$ in the narrow region between $D_2$ and $D_3$, which satisfy that
$$|\nabla u (\mathbf{x}_i)|\geq C\frac{r_1r_3}{r_1+r_3}\frac{1}{\sqrt{r_2}}\frac 1 {\sqrt{\epsilon_i}}~~\mbox{for}~i=1,~2.$$

\par For any entire harmonic function $H$, let $u$ be the solution
to \eqref{eq:conductivity} with $H$. Then, there is a positive
constant $C$ independent of $\epsilon_1$, $\epsilon_2$, $r_1$, $r_2$
and $r_3$ such that
$$|\nabla u (\mathbf{x})|\leq C\frac{r_1r_3}{r_1+r_3}\frac{1}{\sqrt{r_2}}\frac 1 {\sqrt{\epsilon_i}}~~\mbox{for}~i=1,~2.$$
 in the narrow regions between $D_1~\mbox{and}~D_2$, and between $D_2~\mbox{and}~D_3$, respectively.

\end{thm}

In this paper, we first estimate the lower-bounds in terms of the
radii of inclusions. Based on this estimates we derive the
remarkable blow-up rate increasing phenomena when a small conducting
inclusion is located in-between region of two inclusions. This paper
is organized as follows: In section 2, we explain the method to
calculate the potential difference in the case of two disks. We then
derive the lower bound of Case (A) in section 3; Case (B) in section
4. In the case of the upper bounds, the major part of derivation
overlaps in Case (A) and Case (B). Thus, the derivation is presented
in Subsection \ref{upp}. Based on the similar derivation, we can
also obtain the analogues of Theorem \ref{thm:A} and \ref{thm:B} for
the inclusions associated by a sufficiently general class of shapes.

\subsubsection* {Analogues of Theorem \ref{thm:A} and \ref{thm:B} for a sufficiently general class of shapes}
The proofs of Theorem \ref{thm:A} and \ref{thm:B} are flexible
enough even though the results are restricted to circular
inclusions. The estimates presented in Theorem \ref{thm:A} and
\ref{thm:B} can be extended to the inclusions associated by a
sufficiently general class of shapes. To consider a large class of
shapes, we make the geometric assumptions more precise. To define
$D_1$, $D_2$ and $D_3$, we consider three domains
$D_{\mbox{right}}$, $D_{\mbox{center}}$ and $D_{\mbox{left}}$ in $
\mathbb{R}^2$. In addition, we assume that $\varphi_{\mbox{right}} :
\mathbb{C} \backslash B_1(0) \rightarrow \mathbb {R}^2 \backslash
D_{\mbox{right}}$, $\varphi_{\mbox{center}} : \mathbb{C} \backslash
B_1(0) \rightarrow \mathbb {R}^2 \backslash D_{\mbox{center}}$ and
$\varphi_{\mbox{left}} : \mathbb{C} \backslash B_1(0) \rightarrow
\mathbb {R}^2 \backslash D_{\mbox{left}}$ are conformal mappings in
$C^2 ( \mathbb{C} \backslash B_1(0))$ such that
$\varphi_{\mbox{right}}'(z)\neq 0~\mbox{ and}
~\varphi_{\mbox{left}}'(z)\neq 0$ for $z \in
\partial B_1 (0)$. Here, we do not distinguish
$\mathbb{R}^2$ from $\mathbb{C}$.  The $C^2$ regularity condition of
these conformal mappings doses not allow non-smooth inclusions such
as polygons, but Riemann mapping theorem yields a sufficiently
general class of shapes: refer to Ahlfors \cite{A}. Now, we consider
the analogues of Theorem \ref{thm:A} and \ref{thm:B} for two cases
as follows:
\begin{itemize}

 \item[(C)]  One domain and a pair of partially overlapping domains, similarly to Figure \ref{caseCD}: there is a small portion of another domain protruding from a inclusion, i.e.,
$L=2$, and $D_1$ and $D_2$ are $\epsilon$-distanced domains defined
as
\begin{equation}\label{thC} D_1=D_{\mbox{left}}\mbox{ and } D_2= \left ( r_2 D_{\mbox{center}}\right )\cup D_{\mbox{right}},
\end{equation}
where $r_2 D_{\mbox{center}}$ is the $r_2$ times diminished domain
of $D_{\mbox{center}}$. We suppose that $D_2$ is a connected domain,
$\mbox{dist}(D_1,D_2)=\mbox{dist}(D_1,r_2 D_{\mbox{center}})$, $$
\mbox{dist}(D_1,D_{\mbox{right}})\backsimeq r_2,$$ $$D_1 \subset
\mathbb{R_{-}}\times \mathbb{R}~ \mbox{and} ~D_2 \subset
\mathbb{R_{+}}\times \mathbb{R}.$$ In addition, we also assume that
$r_2$ is small enough and $$0< \epsilon \ll r_2,$$ and that the
boundaries $\partial D_1$, $\partial D_2$ and $\partial
D_{\mbox{right}}$ are strictly convex in the narrow region between
$D_1$ and $D_2$.

 \item[(D)]  Three disjoint domains $D_1$, $D_2$ and $D_3$, Figure \ref{caseCD}:
 a small inclusion $D_2$ is disjointly embedded into the in-between area of two other domains, i.e.,
$L=3$, and
\begin{equation}\label{thD} D_1=D_{\mbox{left}},~D_2= r_2 D_{\mbox{center}}~\mbox{ and }
D_3= D_{\mbox{right}} \end{equation}where $r_2 D_{\mbox{center}}$ is
the $r_2$ times diminished domain of $D_{\mbox{center}}$.  We assume
that

 $D_1$ and $D_2$ are $\epsilon_1$ apart, $D_2$ and $D_3$ are $\epsilon_2$
 apart, and $D_1$ distances enough from $D_3$ that $r_2$ is
 sufficiently small and
$$0<\epsilon_{i}\ll r_2~~\mbox{for}~i =1,2,$$
since $D_2$ is regarded as the cross-section of the thin fiber, that
the boundaries $\partial D_1$, $\partial D_2$ and $\partial
D_{\mbox{right}}$ are strictly convex in the narrow region between
$D_1$ and $D_2$,  $$ D_1 \subset \mathbb{R}_- \times \mathbb{R}
~\mbox{and}~ D_2 \cup D_3 \subset \mathbb{R}_+ \times \mathbb{R} .$$

\end{itemize}

\begin{thm}[Case C: Protruding small lump]\label{thm:C}
Let $D_1$ and $D_2$ be defined as \eqref{thC}. Then there is a
positive constant $C$ independent of $\epsilon$ and $r_2$ such that
$$u\Bigr|_{\p D_2}-u\Bigr|_{\p D_1}\geq C\frac{1}{\sqrt{r_2}}\sqrt{\epsilon},$$
where $u$ is the solution to \eqref{eq:conductivity} with
$H(x_1,x_2)=x_1$. As a result, by the Mean Value Theorem, there is a
point $\mathbf{x}_0$ in the narrow region between $D_1$ and $D_2$
such that
$$|\nabla u (\mathbf{x}_0)|\geq C\frac{1}{\sqrt{r_2}}\frac 1 {\sqrt{\epsilon}}.$$

\par For any entire harmonic function $H$, let $u$ be the solution
to \eqref{eq:conductivity} with $H$. Then, there is a positive
constant $C$ independent of  $r_2$ and $\epsilon$ such that
$$|\nabla u (\mathbf{x})|\leq C\frac{1}{\sqrt{r_2}}\frac 1
{\sqrt{\epsilon}}$$ in the narrow region between
$D_1~\mbox{and}~D_2$.
\end{thm}

\begin{thm}[Case D: Disjointly embedded small inclusion] \label{thm:D}
Let $D_i$, $i=1,2,3$, be balls defined as \eqref{thD}. Then there is
a positive constant $C$ independent of $r_2$, $\epsilon_1$ and
$\epsilon_2$ such that
$$u\Bigr|_{\p D_2}-u\Bigr|_{\p D_1}\geq C\frac{1}{\sqrt{r_2}}\sqrt{\epsilon_1},$$
and
$$u\Bigr|_{\p D_3}-u\Bigr|_{\p D_2}\geq C\frac{1}{\sqrt{r_2}}\sqrt{\epsilon_2},$$
where $u$ is the solution to \eqref{eq:conductivity} with
$H(x_1,x_2)=x_1$. As a result, by the Mean Value Theorem, there
exists points; $\mathbf{x}_1$ in the narrow region between $D_1$ and
$D_2$; $\mathbf{x}_2$ in the narrow region between $D_2$ and $D_3$,
which satisfy that
$$|\nabla u (\mathbf{x}_i)|\geq C\frac{1}{\sqrt{r_2}}\frac 1 {\sqrt{\epsilon_i}}~~\mbox{for}~i=1,~2.$$

\par For any entire harmonic function $H$, let $u$ be the solution
to \eqref{eq:conductivity} with $H$. Then, there is a positive
constant $C$ independent of $r_2$,  $\epsilon_1$ and $\epsilon_2$
such that
$$|\nabla u (\mathbf{x}_i)|\leq C\frac{1}{\sqrt{r_2}}\frac 1 {\sqrt{\epsilon_i}}~~\mbox{for}~i=1,~2.$$
 in the narrow regions between
$D_1~\mbox{and}~D_2$, and $D_2~\mbox{and}~D_3$, respectively.
\end{thm}
We derive the lower bound of Case (C) in section 3; Case (D) in
section 4. In the case of the upper bounds, the major part of
derivation overlaps in Case (A), Case (B), Case (C) and Case (D).
Thus, the main idea is presented in Subsection \ref{upp}.

\section{Preliminary}
\subsection{Calculation of the potential difference}\label{section:difference}
We explain the main idea to calculate the difference of potential
between two adjacent, possibly disconnected, conductors.

In this section, differently from \eqref{eq:conductivity}, $D_i$, $i=1,2$, could be also the union of two disjoint domains.
Define $u$ as the solution to \eqref{eq:conductivity}, where it is assigned one constant value throughout $D_i$ even when $D_i$ is disconnected.
Now, define $h$ as the solution to
\begin{equation} \label{def:h}
\quad \left\{
\begin{array}{ll}
\ds\Delta h  = 0,\quad&\mbox{in }{\mathbb{R}^2 \backslash \overline{(D_1 \cup D_2)}}, \\
\ds h= O(|\textbf{x}|^{-1}),\quad&\mbox{as } |\textbf{x}| \rightarrow \infty,\\
\ds h |_{\partial D_i} = k_i\mbox{ (constant)},\quad&\mbox{for }  i = 1, 2,\\
\ds\int_{\partial D_i} {\partial_{  \nu} h }~dS =(-1)^i,\quad&\mbox{for }  i = 1, 2,
\end{array}
\right.
\end{equation}
where $\nu$ is the outward unit normal vector of ${\mathbb{R}^n
\backslash \overline{(D_1 \cup D_2)}}$, i.e., directed inward of $D_i$.
To indicate the dependence of $u$ and $h$ on $D_1$ and
$D_2$, we denote them as
\begin{equation}\label{def:Phi}
u=\Phi[D_1,D_2],
\end{equation}
\begin{equation}\label{def:Psi}
h=\Psi[D_1,D_2].
\end{equation}

The potential difference of $u$ in $D_1$ and $D_2$ is represented in terms of $h$ as follows.
\begin{lem}\label{lemma:h}(\cite{Y})
\begin{equation}\label{eq:102}
u\Bigr|_{\partial D_2} - u\Bigr|_{\partial D_1}=\int_{\partial D_1}
H{\partial_{  \nu} h }\ dS+ \int_{\partial D_2}H {\partial_{ \nu} h
} \ dS.\end{equation}
\end{lem}
The lemma above can be derived by the Divergence Theorem, see
\cite{Y}.

\subsection{Two disks in $\RR^2$}\label{section:twodisks}
Using Lemma \ref{lemma:h}, we can easily calculate the potential difference $u|_{D_2}-u|_{D_1}$ of the solution $u$ to \eqnref{eq:conductivity} when
\begin{equation}\label{twodisks}
D_1 = B_{r_1}(\mathbf{c}_1)~\mbox{and}~D_2 =
B_{r_2}(\mathbf{c}_2),
\end{equation}
where $\mathbf{c}_1=(-r_1 - \frac{\epsilon}{2}, 0)$ and $\mathbf{c}_2=(r_2+ \frac{\epsilon}{2},0)$.

Let $R_i$ be the reflection with respect to $D_i$, in other words,
$$R_i(\mathbf{x})=\frac{r_i^2(\mathbf{x}-\mathbf{c}_i)}{|\mathbf{x}-\mathbf{c}_i|^2}+ \mathbf{c}_i,\ i=1,2,$$
and $\textbf{p}_1\in D_1$ be the fixed point of $R_1 \circ R_2$, then
$R_2(\textbf{p}_1)(=:\textbf{p}_2)$ is the fixed point of  $R_2 \circ R_1$, and
$$\textbf{p}_1 =
\Bigr(-\sqrt{2}\sqrt{\frac{r_1r_2}{r_1+r_2}}\sqrt\epsilon +O(\epsilon)
,0\Bigr)\mbox{ and }\textbf{p}_2 =
\Bigr(\sqrt{2}\sqrt{\frac{r_1r_2}{r_1+r_2}}\sqrt\epsilon +O(\epsilon)
,0\Bigr).$$
Moreover, we can easily show that
\begin{equation}\label{h:twodisks}\Psi[D_1,D_2]=\frac 1 {2\pi} \left( \log |\textbf{x}-\textbf{p}_1|- \log |\textbf{x}-\textbf{p}_2|\right).\end{equation}
By an elementary calculation, it can be shown that the middle point
$\frac {\textbf{p}_1 + \textbf{p}_2}{2}$ exists between two
approaching points $\left(-\frac {\epsilon} 2, 0 \right)$ and
$\left(\frac {\epsilon} 2, 0 \right)$. Applying the middle point
property to estimate for $\Psi[D_1,D_2]\left(\pm\frac {\epsilon} 2,
0 \right)$, we can get the following lemma.
\begin{lem}\label{lemma:twodisks} There is a constant $C >0$
independent of $\epsilon$, $r_1$ and $r_2$ such that
$$\frac 1 C \sqrt{\frac {r_1+r_2} {r_1r_2}}\sqrt\epsilon \leq \Psi[D_1,D_2]\big|_{\partial D_2} - \Psi[D_1,D_2]\big|_{\partial D_1} \leq  C \sqrt{\frac {r_1+r_2} {r_1r_2}}\sqrt\epsilon $$
for small $\epsilon >0$. \end{lem} From Lemma \ref{lemma:h}, we
calculate the potential difference of $u$.
\begin{lem} \label{lem++}
Let $H(x_1,x_2)$ be an entire harmonic function. The solution $u$ to
\eqnref{eq:conductivity} where $L=2$ and $D_l$, $l=1,2$, are given
as \eqref{twodisks}
 satisfies
 \be
 \begin{array}{ll}u|_{\partial D_2} - u|_{\partial D_1} &=
H(\textbf{p}_2)-H(-\textbf{p}_1)\\ &= 2\sqrt{2} \partial_{x_1} H (0,0)
\sqrt{\frac{r_1r_2}{r_1+r_2}}\sqrt\epsilon+O(\epsilon).
 \label{eq:107}\end{array}\ee

 \end{lem}
 \begin{rem}
 Referring to the mean value theorem, there exists a point
$\textbf{x}_2$ between $\partial D_1$ and $\partial D_2$ such that
\be |\nabla u (\textbf{x}_2)|\geq2\sqrt{2} |\partial_{x_1} H(0,0)|
\sqrt{\frac{r_1r_2}{r_1+r_2}} \frac 1 {\sqrt{\epsilon}}.
\label{eq:108}\ee for any sufficiently small $\epsilon >0$.
Moreover, as a result in \cite{AKLLL}, there is a constant $C$
independent of $\epsilon$, $r_1$ and $r_2$ such that
$$ \norm{\nabla u }_{L^{\infty}(\Omega \setminus (D_1 \cup D_2)) }\leq C \norm {\nabla H}_{L^{\infty}(\Omega)}
\sqrt{\frac{r_1r_2}{r_1+r_2}} \frac 1 {\sqrt{\epsilon}} $$ where
$\Omega=B_{4(r_1+r_2)}(0,0).$
\end{rem}

\section{One disk and a pair of partially overlapping disks}
In this section, we consider two $\epsilon$-distanced domains $D_1$
and $D_2$, see Case (A) at Figure \ref{caseAB}, where
\begin{equation*}D_1=B_{r_1}(\mathbf{c}_1)\mbox{ and }D_2=B_{r_2}(\mathbf{c}_2)\cup B_{r_3}(\mathbf{c}_3),
\end{equation*}
where \begin{equation*}\mathbf{c}_1=(-r_1-\frac{\epsilon}{2},0),\ \mathbf{c}_2=(r_2+\frac{\epsilon}{2},0),\mbox{ and } \mathbf{c}_3=(r_3+a+\frac{\epsilon}{2},0).\end{equation*}
Here, $B_{r_2}(\mathbf{c}_2)$ is a small lump of $B_{r_3}(\mathbf{c}_3)$, and we assume
$$
B_{r_2}(\mathbf{c}_2)\cap B_{r_3}(\mathbf{c}_3)\neq\emptyset,\mbox{ i.e., }0<a<2r_2,$$
and $$\mbox{ and }0<\epsilon\ll r_2\ll \min(r_1,\ r_3).$$

Define
\begin{equation}\label{h:D1D2}h=\Psi[D_1, D_2],\mbox{ and } u=\Phi[D_1, D_2],
\end{equation}
and
 \begin{equation}\label{hj:D1D2}h_j=\Psi[D_1, B_{r_j}(\mathbf{c}_j)],\mbox{ and }u_j=\Phi[D_1, B_{r_j}(\mathbf{c}_j)],\qquad j=2,3,
\end{equation}
where $\Psi$ and $\Phi$ defined in section \ref{section:difference}.

\subsection{Properties of $h$ and $h_j$, $j=2,3$}

\begin{lem}\label{property:hi}
Let $h=\Psi[D_1, D_2]$, then we have
\begin{equation}\label{nuh:small}
\p_\nu h(x)=O(\sqrt{\epsilon}),\quad x\in \p D_2\setminus B_{r_2}(\mathbf{c}_2).
\end{equation}
\end{lem}
\pf
Set $$B_i=B_{r_i}(\mathbf{c}_i),\quad\mbox{for }i=1,2,3.$$

We choose a smooth domain $\widetilde{\Omega}$ as
follows: \begin{align*}& \widetilde{\Omega} \subset (B_2 \cup
B_3),~B_3 \subset \widetilde{\Omega}\notag\\&\partial
\widetilde{\Omega} \setminus B_2=\partial B_3 \setminus B_2\notag\\
& (\partial \widetilde{\Omega} \cap \partial B_2)\setminus
\partial B_3 = (\frac 1 2\epsilon,0)\notag\end{align*}
Let $\widetilde{h} = \Psi[B_1, \widetilde{\Omega}]$. Then, we
consider $V$ defined in $\mathbb{R}^2 \setminus (B_1 \cup B_2 \cup
B_3)$ as follows:
$$V = h - \frac {h|_{\partial B_1} - h|_{\partial (B_2 \cup B_3)}}{\widetilde{h}|_{\partial B_1} - \widetilde{h}|_{\partial \widetilde{\Omega}}} \widetilde{h}.$$
Then, it follows that
$$V|_{\partial B_1} = V|_{\partial \widetilde{\Omega}\setminus B_2}=\mbox{a constant} .$$
Since $h|_{\partial (B_2 \cup B_3)}>h|_{\partial B_1}$ and $
\widetilde{h}|_{\partial
\widetilde{\Omega}}>\widetilde{h}|_{\partial B_1}$, the minimum of
$V$ attains on $\partial \widetilde{\Omega} \setminus B_2$ and
$\partial B_1$. Thus, we have \be\partial_{\nu} h - \frac
{h|_{\partial B_1} - h|_{\partial (B_2 \cup
B_3)}}{\widetilde{h}|_{\partial B_1} - \widetilde{h}|_{\partial
\widetilde{\Omega}}} \partial_{\nu} \widetilde{h}\leq
0~\mbox{on}~\partial B_1 \cup (\partial \widetilde{\Omega}\setminus
B_2).\label{eq:103}\ee By the integration on $\partial B_1$, we have
$$0<\frac
{h|_{\partial B_1} - h|_{\partial (B_2 \cup
B_3)}}{\widetilde{h}|_{\partial B_1} - \widetilde{h}|_{\partial
\widetilde{\Omega}}}\leq 1.$$ Using the bound \eqref{eq:103} once
more, we have
$$\partial_{\nu} h  \leq  \partial_{\nu}\widetilde{h}|_{\partial
\widetilde{\Omega}} ~\mbox{on}~\partial \widetilde{\Omega}\setminus
B_2=\partial B_3 \setminus B_2.$$

The domain $\widetilde{ \Omega}$ is smooth so that we can use the
method presented by Yun \cite{Y,Y2}. Then, up to a conformal mapping
to a circle, $\partial_{\nu} h$ is bounded by constant times the
Poisson Kernel with respect to a interior point $\sqrt \epsilon$
distanced from the boundary (refer to the inequality (9) in
\cite{Y2}). Note that $\partial B_3 \setminus B_2$ distances enough
from $(\epsilon,0)$. Thus, we have
$$\partial_{\nu} h  \leq  \partial_{\nu}\widetilde{h}|_{\partial
\widetilde{\Omega}} \leq C \sqrt {\epsilon} ~\mbox{on}~\partial B_3
\setminus B_2.$$ Therefore, we have completed the proof of the
lemma. \qed

\begin{lem}\label{lem:Anuh}
\begin{equation}
\p_\nu h (x)\leq M\p_\nu h_3(x),\quad x\in\p D_1,
\end{equation}
where\begin{equation}\label{def:M}M=\frac{h\bigr|_{\p
D_2}-h\bigr|_{\p D_1}}{h_3\bigr|_{\p
B_{r_3}(\mathbf{c}_3)}-h_3\bigr|_{\p D_1}}.\end{equation}
\end{lem}
\pf
 Define $$W=h-Mh_3,\qquad \mbox{in }\RR^2\setminus(D_1\cup D_2).$$
Since $h$ is constant on $\p  D_2$, $M>0$, and $h_3$ takes it's
maximum on $\p B_{r_3}(\mathbf{c}_3)$,
$$W\Bigr|_{\p D_2\setminus B_{r_2}(\mathbf{c}_2)}-W\Bigr|_{\p D_2\setminus B_{r_3}(\mathbf{c}_3)}=-M(h_3\Bigr|_{\p B_{r_3}(\mathbf{c}_3)}-h_3\Bigr|_{\p D_2\setminus B_{r_3}(\mathbf{c}_3)})<0,$$
and $$W\Bigr|_{\p D_2\setminus B_{r_2}(\mathbf{c}_2)}-W\Bigr|_{\p
D_1}=h\Bigr|_{\p D_2}-h\Bigr|_{\p D_1}-M(h_3\Bigr|_{\p
B_{r_3}(\mathbf{c}_3)}-h_3\Bigr|_{\p D_1})=0.$$

Therefore, $W$ takes its minimum on $\p D_1$, and
$$\p_\nu W\leq0,\quad \mbox{on }\p D_1.$$
\qed

\begin{lem}\label{lem:Asamediffer}
\begin{equation}\label{samediffer}
\ds h\Bigr|_{\p D_2}-h\Bigr|_{\p  D_1}=h_2\Bigr|_{\p B_{r_2}(\mathbf{c}_2)}-h_2\Bigr|_{\p D_1}+O(\epsilon).
\end{equation}
\end{lem}
\pf
 Note that $$\int_{\p D_1}\p_\nu (h-h_2)\ dS=0,$$ and
 \begin{align*}\ds\int_{\p D_2}\p_\nu(h-h_2)\ dS&=\ds\int_{\p D_2}\p_\nu h\ dS-\int_{\p B_{r_2}(\mathbf{c}_2)}\p_\nu h_2\ dS-\int_{\p (D_2\setminus B_{r_2}(\mathbf{c}_2))}\p_\nu h_2 \ dS\\&=\ds 1-1-0=0.\end{align*}
With the fact that $h|_{\p D_1}$ and $h|_{\p D_2}$ are constants and
the (exterior) Divergence Theorem, we have that
\begin{align*}
\ds0&=\int_{\p D_1}\p_\nu(h-h_2)h\ dS+\int_{\p D_2}\p_\nu (h-h_2)h\ dS\\
\ds&=\int_{\p D_1}(h-h_2)\p_\nu h\ dS+\int_{\p D_2}(h-h_2)\p_\nu h\ dS.
\end{align*}
Hence,
\begin{align*}
h\Bigr|_{\p D_1}-h\Bigr|_{\p  D_2}
&=\int_{\p D_1}h\p_\nu h\ d S+\int_{\p D_2}h\p_\nu h\ d S\\
&=\int_{\p D_1}h_2\p_\nu h \ dS+\int_{\p D_2}h_2\p_\nu h\ dS\\
&=h_2\Bigr|_{\p D_1}-h_2\Bigr|_{\p B_{r_2}(\mathbf{c}_2)}+\int_{\p D_2}(h_2-h_2\Bigr|_{\p B_{r_2}(\mathbf{c}_2)})\p_\nu h\ dS
\end{align*}
From \eqref{h:twodisks}, there is a constant $C$ dependent of $a$, see \eqref{eqn:lumpcenters}, such that
$$\Bigr|(h_2-h_2\Bigr|_{\p B_{r_2}(\mathbf{c}_2)})(x)\Bigr|\leq C\sqrt{\epsilon},\quad\mbox{for all }x\in \p  D_2\setminus B_{r_2}(\mathbf{c}_2).$$
Therefore, with \eqref{nuh:small} as well, we obtain \eqref{samediffer}.
\qed

\subsection{Proof Theorem \ref{thm:A}}
Let $H(x_1,x_2)=x_1$ and  $\nu$ be the unit normal vector of ${\mathbb{R}^2\backslash \overline{(D_1 \cup D_2)}}$, i.e., directed inward to $D_i$, $i=1,2$. Remind that we defined
\begin{equation}\label{h:D1D2}h=\Psi[D_1, D_2],\mbox{ and } u=\Phi[D_1, D_2],
\end{equation}
and
 \begin{equation}\label{hj:D1D2}h_j=\Psi[D_1, B_{r_j}(\mathbf{c}_j)],\mbox{ and }u_j=\Phi[D_1, B_{r_j}(\mathbf{c}_j)],\qquad j=2,3,
\end{equation}
where $\Psi$ and $\Phi$ defined in section \ref{section:difference}.

Note that $\partial_{ \nu}h\Bigr|_{\p D_2}<0$, $H <0$ on $\p D_1$
and $H >0$ on $\p D_2$, and, as a result, from Lemma \ref{lemma:h},
we have
\begin{align}
u\Bigr|_{\partial D_2} - u\Bigr|_{\partial D_1}&=\int_{\partial D_1}\nonumber
({\partial_{  \nu} h }) H\ dS+ \int_{\partial D_2}( {\partial_{ \nu}
h }) H\ dS\\
&\geq\int_{\p D_1}H\p_\nu h\ dS.\label{ud1d2_first}
\end{align}

Applying the lemma \ref{lem:Anuh}, Lemma \ref{lem:Asamediffer}, \eqref{ud1d2_first} becomes \begin{align*}
\ds u\Bigr|_{\p D_2}-u\Bigr|_{\p D_1}&\geq \frac{h\bigr|_{\p D_2}-h\bigr|_{\p D_1}}{h_3\bigr|_{\p B_{r_3}(\mathbf{c}_3)}-h_3\bigr|_{\p D_1}}\int_{\p D_1} H\p_\nu h_3\ dS\\
\ds&\geq \frac{h_2\bigr|_{\p B_{r_2}(\mathbf{c}_2)}-h_2\bigr|_{\p D_1}+O(\epsilon)}{h_3\bigr|_{\p B_{r_3}(\mathbf{c}_3)}-h_3\bigr|_{\p D_1}}\sqrt{2}\sqrt{\frac{r_1r_3}{r_1+r_3}r_2}.
\end{align*}
It follows  from Lemma \ref{lemma:twodisks} that
$$h_2\Bigr|_{\p B_{r_2}(\mathbf{c}_2)}-h_2\Bigr|_{\p D_1} \geq C \sqrt{\frac{r_1+r_2} {r_1r_2} } \sqrt{\epsilon}+O(\epsilon)$$
and
$$h_3\Bigr|_{\p B_{r_3}(\mathbf{c}_3)}-h_3\Bigr|_{\p D_1} \leq C\sqrt{\frac{r_1+r_3} {r_1r_3} } \sqrt{r_2}+O(r_2).$$
Therefore,
\begin{align*}
\ds u\Bigr|_{\p D_2}-u\Bigr|_{\p D_1}&\geq C\frac{r_1r_3}{r_1+r_3}\frac{1}{\sqrt{r_2}}\sqrt{\epsilon}.
\end{align*}
This proves  Theorem \ref{thm:A}.
\qed

\subsection{Proof Theorem \ref{thm:C}}
We consider the general shaped domain in Theorem \ref{thm:C}. But,
we take an advantage of the properties of circular inclusions. To
make a connection between circular domains and general shaped
domains, we need to establish the monotonic property of $\Psi$ as
follows:

\begin{lem} \label{lem:mono} [Monotonic property of $\Psi$] Let $D_A$, $D_B$, $\widetilde{D}_A$ and
$\widetilde{D}_B$ be  domains. Assume that $$D_A \subseteq
\widetilde{D}_A~\mbox{and}~D_B \subseteq \widetilde{D}_B. $$ Then,
we have
$$ 0 \leq \Psi [\widetilde{D}_A, \widetilde{D}_B]\big|_{\partial \widetilde{D}_B} -  \Psi [\widetilde{D}_A, \widetilde{D}_B]\big|_{\partial \widetilde{D}_A}
 \leq \Psi [D_A, D_B]\big|_{\partial D_B} -  \Psi [D_A, D_B]\big|_{\partial D_A}.$$
\end{lem}
\pf Without any loss of generality, we consider only the case of
$D_A = \widetilde{D}_A.$ Let $$G =\Psi [{D}_A, \widetilde{D}_B] - M
\Psi [{D}_A, {D}_B] $$ where
$$M = \frac {\Psi [D_A, \widetilde{D}_B]\big|_{\partial D_B} -  \Psi [D_A, \widetilde{D}_B]\big|_{\partial D_A} }{\Psi [{D}_A, {D}_B]\big|_{\partial {D}_B} -  \Psi [{D}_A,{ D}_B]\big|_{\partial {D}_A} } .$$
The minimum  of $G$ attains on $\partial D_A$. By the Hopf's Lemma,
we have
$$\partial_{\nu} G \leq 0~\mbox{on}~\partial D_A.$$
Integrating $\p_{\nu} G$ on $\p D_A$, we have $-1 + M \leq 0$.
Therefore, we have
$$ 0 \leq \Psi [{D}_A, \widetilde{D}_B]\big|_{\partial \widetilde{D}_B} -  \Psi [{D}_A, \widetilde{D}_B]\big|_{\partial {D}_A}
 \leq \Psi [D_A, D_B]\big|_{\partial D_B} -  \Psi [D_A, D_B]\big|_{\partial D_A}.$$
Repeating the same argument again, we can obtain the disable
inequality.

\qed

Applying $D_{\mbox{left}}$, $r_2 D_{\mbox{center}}$ and
$D_{\mbox{right}}$ instead of $B_{r_i}(\mathbf{c}_i)$, $i=1,2,3$, to
the argument presented in the proof of Theorem \ref{thm:A}, we can
obtain
$$
\ds u\Bigr|_{\p r_2 D_{\mbox{center}}}-u\Bigr|_{\p
D_{\mbox{left}}}\geq C \frac{h_2\bigr|_{\p\left( r_2
D_{\mbox{center}} \right) }-h_2\bigr|_{\p
D_{\mbox{left}}}+O(\epsilon)}{h_3\bigr|_{\p
D_{\mbox{right}}}-h_3\bigr|_{\p D_{\mbox{left}}}}\sqrt{r_2} $$when
$H(x_1,x_2) = x_1$. Here, $h_1 = \Psi [D_{\mbox{left}}, r_2
D_{\mbox{center}}]$ and $h_2 = \Psi [D_{\mbox{left}},
D_{\mbox{right}}]$.

 It follows that from Yun
\cite{Y, Y2} that
$$h_3\Bigr|_{\p D_{\mbox{right}}}-h_3\Bigr|_{\p D_{\mbox{left}}}\simeq  \sqrt{r_2}.$$
To estimate $h_2\Bigr|_{\p r_2 D_{\mbox{center}}}-h_2\Bigr|_{\p
D_{\mbox{left}}}$, we choose two disks $B_{\mbox{left}}$ and
$B_{\mbox{center}}$ containing $D_{\mbox{left}}$ and
$D_{\mbox{center}}$ such that the distance between $B_{\mbox{left}}$
and $r_2 B_{\mbox{center}}$ is $\epsilon$. Using Lemma
\ref{lem:mono} and \ref{lemma:twodisks}, we have
$$h_2\Bigr|_{\p \left(r_2 D_{\mbox{center}}\right)}-h_2\Bigr|_{\p D_{\mbox{left}}} \gtrsim  \sqrt{\frac {\epsilon}{r_2}}.$$
Note that $D_1= D_{\mbox{left}}$, $D_2 = r_2 D_{\mbox{center}}$ and
$D_3 = D_{\mbox{right}}$ in this theorem.  Therefore,
\begin{align*}
\ds u\Bigr|_{\p D_2}-u\Bigr|_{\p D_1}&\geq C
\frac{1}{\sqrt{r_2}}\sqrt{\epsilon}.
\end{align*}
This proves  Theorem \ref{thm:C}. \qed

\smallskip

\section{Three disjoint smooth domains}\label{section:threedomains}
We consider three disjoint inclusion case, see Figure \ref{caseAB}
and \ref{caseCD}, a small one is disjointly embedded into the
in-between area of two others, and prove Theorem \ref{thm:B} and
\ref{thm:D}. We assume that $D_1$ and $D_2$ are closely spaced with
the distance $\epsilon_1$, and $D_2$ and $D_3$ are closely space
with $\epsilon_2$, but $D_1$ and $D_3$ are not close, and that
$D_1$, $D_2$ and $D_3$ have the boundary regularity given in Theorem
\ref{thm:D}.

\subsection{Solution representation of $u$}\label{subsec}
Let $H^c$ be a harmonic function outside of $\cup_{i=1}^3 D_i$ and
have the same constant value in $\cup_{i=1}^3 D_i$ satisfying that
\begin{equation} \label{def:Hc}
\quad \left\{
\begin{array}{ll}
\ds\Delta H^c  = 0,\quad&\mbox{in }{\mathbb{R}^2 \backslash \overline{\cup_{i=1}^3 D_i}},\\
\ds H^c(\textbf{x})- H(\textbf{x}) = O(|\textbf{x}|^{-1}),\quad&\mbox{as } |\textbf{x}| \rightarrow \infty,\\
\ds H^c |_{\cup_{i=1}^3\partial D_i} = C_H~ \mbox{(constant)}.
\end{array}
\right.
\end{equation}
Since $H^c-H$ is harmonic at infinity, $H^c-H$ attains maximum only at the boundary points of $D_i$, $i=1,2,3$.
To make $H^c-H$ attains zero at infinity, $C_H$ should satisfy \begin{equation}\label{range:c}
\ds-\bigr\|H\bigr\|_{L^\infty{(\cup_{i=1}^3 D_i})}\leq C_H\leq\bigr\|H\bigr\|_{L^\infty{(\cup_{i=1}^3 D_i})}.
\end{equation}
Moreover, $H^c$ satisfies $\sum_{i=1}^3\int_{\partial D_i} {\partial_{  \nu} H^c }~dS = 0.$

The solution $u$ to \eqref{eq:conductivity} is represented as
\begin{equation}\label{rep:u3balls}
\ds u(\textbf{x})=H^c(\textbf{x})+c_1h_1(\textbf{x})+c_2h_2(\textbf{x}),
\end{equation}
where
$$h_1=\Psi\bigr[D_1,(D_2\cup D_3)\bigr],\ h_2=\Psi\bigr[(D_1\cup D_2),D_3\bigr],$$
and
\begin{equation}\label{def:c}
\left(\begin{array}{c}
\ds c_1\\
\ds c_2
\end{array}\right)
=\ds-\left( \begin{array}{cc}
\ds -1 &\ds\int_{\p D_1}\p_\nu h_2\ d S\\
\ds \int_{\p D_2}\p_\nu h_1\ dS&\ds\int_{\p D_2}\p_\nu h_2\ d S \end{array} \right)^{-1}
\left(\begin{array}{c}
\ds\int_{\p D_1}\p_\nu H^c\ dS\\
\ds\int_{\p D_2}\p_\nu H^c\ dS
\end{array}\right),
\end{equation}
where $\Psi$ is defined as \eqref{def:h} and \eqref{def:Psi}. The
equality \eqref{def:c} is from the integration of $\partial_{\nu } u
$ on $\p D_1$ and  $\p D_2$.

 Applying the upper bound on the gradient of  solution without the potential difference among the boundaries to conductivity equation derived in Bao et al.
\cite{BLY}, we can show that $\nabla H^c$ does not blow-up (also
refer to \cite{Y}). Using Lemma \ref {lem:B} in the following
section, we have $\int_{\p D_2}\p_\nu h_1\ dS = 1+  O(\sqrt
{\epsilon_1})$. This implies that
\begin{equation*}
\left(\begin{array}{c}
\ds c_1\\
\ds c_2
\end{array}\right)
\thickapprox\ds-\left( \begin{array}{cc}
\ds -1 &\ds 0 \\
\ds 1&\ds-1
\end{array} \right)^{-1} \left(\begin{array}{c}
\ds\int_{\p D_1}\p_\nu H^c\ dS\\
\ds\int_{\p D_2}\p_\nu H^c\ dS
\end{array}\right).
\end{equation*}
Thus, the coefficient $c_i$, $i=1,2$, is bounded independently of
$\epsilon_1$ and $\epsilon_2$. Therefore, the blow-up rate of
$\nabla u$ essentially relies on $\nabla h_i$. In this respect, we
consider the properties of $h_i$ in the following section.

\subsection{Properties of $h_1$ and $h_2$}
We build the optimal bounds of $u$ based on \eqref{rep:u3balls}; it
is essential to drive properties of $h_1$ and $h_2$ in the narrow
regions between inclusions. Let $h_1$ and $h_2$ be as follows:
$$h_1 = \Psi\bigr[D_1,(D_2\cup
D_3)\bigr]~\mbox{and}~h_2=\Psi\bigr[(D_1\cup D_2),D_3\bigr]. $$

\begin{prop}\label{prop:nablabound} There are the
following  estimates for $h_1$ and $h_2$:
\begin{itemize}
\item[(i)]  In the narrow region between $D_1$ and $D_2$, we have
$$\nabla h_1=O\Bigr(\frac{1}{\sqrt{\epsilon_1}}\Bigr)~\mbox{and}~\nabla h_2=O(\sqrt{\epsilon_2}).$$
\item[(ii)] In the narrow region between $D_2$ and $D_3$, we have
$$\nabla h_1=O(\sqrt{\epsilon_1})~\mbox{and}~\nabla h_2=O\Bigr(\frac{1}{\sqrt{\epsilon_2}}\Bigr).$$
\item[(iii)] $$h_1 |_{\partial D_2 \cup \partial D_3 } - h_1 |_{\partial D_1} \simeq \sqrt {\epsilon_1}$$ and $$  h_2 |_{\partial D_3} - h_2 |_{\partial D_1 \cup \partial D_2 } \simeq \sqrt {\epsilon_2}.$$
\end{itemize}
\end{prop}

\pf We consider $\nabla h_1$. By Lemma \ref{lem:A} and \ref{lem:C},
we have

$$0>
\partial_{\nu} h_1 \geq C\p_{\nu}\Psi[D_1, D_4]~\mbox{on}~\partial D_1$$
and
$$0 <
\partial_{\nu} h_1 \leq \partial_{\nu}
\Psi[D_1,D_2]~\mbox{on}~\partial D_2,$$ and by Lemma \ref{lem:B},
$$0 <
|\partial_{\nu} h_1| \leq C \sqrt {\epsilon_1}~\mbox{on}~\partial
D_3,$$ where $D_4$ is defined in Lemma \ref{lem:C}.  Without any
loss of generality, we assume that
$$\left(-\frac {\epsilon_1} 2, 0\right) \in \p D_1,~\left(\frac {\epsilon_1} 2, 0\right) \in \p D_2~\mbox{and}~\mbox{dist}(D_1, D_2) =\epsilon_1.$$
Let
$$p(\mathbf{x})= \log |\mathbf{x} - (\sqrt {\epsilon_1},0)| - \log |\mathbf{x} + (\sqrt {\epsilon_1},0) |. $$
Referring to the inequality (9) in \cite{Y2}, there is a constant
$C_1$ such that
$$0< |\nabla h_1| \leq C_1 |\nabla p|~\mbox{on}~\partial (D_1 \cup D_2 \cup D_3).$$
Regarding $(x_1,x_2)$ as a complex number $z=x_1 + x_2 i$, we
consider
$$\rho (z) = \frac {\p_1 h_1(z) - \p_2 h_1(z) i }{C_1 (\p_1 p(z) - \p_2 p(z) i) }.$$
Then, $\rho(z)$ can be extended to $\infty$ as an analytic function.
From definition, $|\rho(z)| < 1 $ on $\partial D_1 \cup \partial D_2
\cup\partial D_3 $. By the maximum principle,  $$|\rho(z)| < 1
~\mbox{in}~  \mathbb{C} \setminus (D_1 \cup  D_2 \cup D_3).$$ Thus,
we have
$$|\nabla h | \leq C_1 |\nabla p | ~\mbox{in}~\mathbb{R}^2 \setminus (D_1 \cup  D_2 \cup D_3). $$
Therefore, $\nabla h_1=O\Bigr(\frac{1}{\sqrt{\epsilon_1}}\Bigr)$ in
the narrow region between $D_1$ and $D_2$, and $\nabla
h_1=O(\sqrt{\epsilon_1})$ in the narrow region between $D_2$ and
$D_3$. Similarly, we have $\nabla
h_2=O\Bigr(\frac{1}{\sqrt{\epsilon_2}}\Bigr)$ in the narrow region
between $D_2$ and $D_3$, and $\nabla h_2=O(\sqrt{\epsilon_2})$ in
the narrow region between $D_1$ and $D_2$. We have proven (i) and
(ii).
\par The estimate (iii) is presented by Lemma \ref{lem:C}.

\qed

\begin{lem}\label{lem:A} We have the following properties:
\begin{itemize}
\item[(i)]$$0< h_1 |_{\partial D_2} - h_1 |_{\partial D_1} \leq \Psi[D_1, D_2]\Big|_{\partial D_2} - \Psi[D_1, D_2]\Big|_{\partial D_1}.$$
\item[(ii)] $$0 < \partial_{\nu} h_1 \leq \partial_{\nu} \Psi[D_1,D_2]~\mbox{on}~\partial D_2.$$
\end{itemize}
\end{lem}
\pf Let $$M=\frac{\Psi[D_1,(D_2\cup D_3)]\Bigr|_{\p
D_1}-\Psi[D_1,(D_2\cup D_3)]\Bigr|_{\p D_2\cup\p
D_3}}{\Psi[D_1,D_2]\Bigr|_{\p D_1}-\Psi[D_1,D_2]\Bigr|_{\p D_2}},$$
and
\begin{align*}\ds G(\mathbf{x})&=\Psi[D_1,(D_2\cup D_3)](\mathbf{x})-\Psi[D_1,(D_2\cup D_3)]\Bigr|_{\p D_2\cup \p D_3}\\
&\qquad\qquad\ds-M\left(\Psi[D_1,D_2](\mathbf{x})-\Psi[D_1,D_2]\Bigr|_{\p
D_2}\right).\end{align*} Then, $G= 0$ on $\p D_1\cup \p D_2$, and
$G>0$ on $\p D_3$. By Hopf's lemma, $$\p_\nu G<0\quad\mbox{on }\p
D_1.$$ This means that \be\partial_{\nu} h_1 \leq M
\partial_{\nu} \Psi[D_1, D_2]~\mbox{on}~\partial D_1.\label{eq:3}\ee
 Note that $h_1 = \Psi[D_1,(D_2\cup
D_3)]$. By integrating $G$ on $\p D_1$, we have the inequality (i).

\par On the other hand, by Hopf's lemma, $$\p_\nu G<0\quad\mbox{on
}\p D_2.$$  This means that
$$\partial_{\nu} h_1 \leq M \partial_{\nu} \Psi[D_1, D_2]~\mbox{on}~\partial D_2.$$
From the inequality (i), $M<1$. Therefore, we have (ii).

\qed

\begin{lem} \label{lem:B} There is a constant $C$ such that
$$0\leq \partial_{\nu} h_1 \leq C \sqrt \epsilon_1 ~\mbox{on}~\partial D_3.$$
\end{lem}
\pf We use the method similar to Lemma \ref{lem:A}. Let
$$M=\frac{\Psi[D_1,(D_2\cup D_3)]\Bigr|_{\p D_1}-\Psi[D_1,(D_2\cup
D_3)]\Bigr|_{\p D_2\cup\p D_3}}{\Psi[D_1,D_3]\Bigr|_{\p
D_1}-\Psi[D_1,D_3]\Bigr|_{\p D_3}},$$ and
\begin{align*}\ds G(\mathbf{x})&=\Psi[D_1,(D_2\cup D_3)](\mathbf{x})-\Psi[D_1,(D_2\cup D_3)]\Bigr|_{\p D_2\cup \p D_3}\\
&\qquad\qquad\ds-M\left(\Psi[D_1,D_3](\mathbf{x})-\Psi[D_1,D_3]\Bigr|_{\p
D_3}\right).\end{align*} Then, $G= 0$ on $\p D_1\cup \p D_3$, and
$G>0$ on $\p D_2$. By Hopf's lemma, $$\p_\nu G<0\quad\mbox{on }\p
D_3.$$ Since $h_1 = \Psi[D_1,(D_2\cup D_3)]$, this inequality means
that
$$0\leq \partial_{\nu} h_1 \leq M \p_{\nu} \Psi[D_1,D_3]~\mbox{on}~\p D_3.$$

\par Now, we estimate the gradient of $M \Psi[D_1,D_3]$. To do so,
we consider the potential difference between $\p D_1$ and $\p D_3$
as follows:
\begin{align*}
M \Psi[D_1,D_3]\Big|_{\partial D_3} -M \Psi[D_1,D_3]\Big|_{\partial
D_1} &= h |_{\partial D_3} -
h|_{\partial D_1}\\&=h |_{\partial D_2} - h|_{\partial D_1}\\
&\leq \Psi[D_1,D_2]\Big|_{\partial D_2} -
\Psi[D_1,D_2]\Big|_{\partial D_1}\\
&\leq C \sqrt {\epsilon_1}\end{align*} The last inequality above was
proven by Yun in his paper \cite{Y,Y2}, since $\Psi[D_1,D_2]$ is
only for two domains. Note that $D_3$ is not close to $D_2$. Owing
to the method in Bao et al. \cite{BLY}, we have
$$\norm {\p_{\nu} M \Psi[D_1,D_3] }_{L^{\infty} (\partial D_3)} \leq C \sqrt {\epsilon_1}.$$
Therefore, we can obtain the result.
 \qed

\begin{lem} \label{lem:C} Let $D_4$ is a disk containing $D_2$ and
$D_3$ with $$\mbox{dist}(D_1,D_4) = \mbox{dist}(D_1,D_2).$$
\begin{itemize}
\item[(i)]There is a positive constant $C$ such that $$0> \partial_{\nu} h_1 \geq C\partial_{\nu} \Psi[D_1, D_4]~\mbox{on}~\partial D_1.$$
\item[(ii)] $$h_1 |_{\partial D_2} - h_1 |_{\partial D_1} \geq \Psi[D_1, D_4]\Big|_{\partial D_4} - \Psi[D_1, D_4]\Big|_{\partial D_1}.$$

\item[(iii)] $$h_1 |_{\partial D_2 \cup \partial D_3 } - h_1 |_{\partial D_1} \simeq \sqrt {\epsilon_1}.$$
\end{itemize}
\end{lem}
\pf To prove (i) and (ii), we use the same derivation to Lemma
\ref{lem:A}. So, we set $$M=\frac{\Psi[D_1,(D_2\cup D_3)]\Bigr|_{\p
D_1}-\Psi[D_1,(D_2\cup D_3)]\Bigr|_{\p D_2\cup\p
D_3}}{\Psi[D_1,D_4]\Bigr|_{\p D_1}-\Psi[D_1,D_4]\Bigr|_{\p D_4}},$$
and
\begin{align*}\ds G(\mathbf{x})&=\Psi[D_1,(D_2\cup D_3)](\mathbf{x})-\Psi[D_1,(D_2\cup D_3)]\Bigr|_{\p D_2\cup \p D_3}\\
&\qquad\qquad\ds-M\left(\Psi[D_1,D_4](\mathbf{x})-\Psi[D_1,D_4]\Bigr|_{\p
D_4}\right).\end{align*} Then $G \Big|_{\p D_1}= 0$ and $G \leq 0$
on ${\p D_4}$. By Hopf's lemma, we have
$$\p_{\nu}G >0~\mbox{on}~\p D_1.$$
By the integration on $\p D_1$, we have (ii) and $M<1$. Therefore,
the inequality $\p_{\nu}G >0$ can also yield (i).

\par From (i) of Lemma \ref{lem:A} and (ii) in this lemma, we have
$$ h_1 |_{\partial D_2} - h_1 |_{\partial D_1} \geq \Psi[D_1,
D_4]\Big|_{\partial D_4} - \Psi[D_1, D_4]\Big|_{\partial D_1}$$ and
$$ h_1 |_{\partial D_2} - h_1 |_{\partial D_1} \leq \Psi[D_1,
D_2]\Big|_{\partial D_2} - \Psi[D_1, D_2]\Big|_{\partial D_1}.$$ The
potential $\Psi[D_1, D_i] ~(i=1,4)$ is only for two domains and
thus, its difference between $D_1$ and $D_i$ ($i=1,2$) was already
estimated in Yun \cite{Y, Y2} as follows: for $i=1,2$,
$$\Psi[D_1,
D_i]\Big|_{\partial D_i} - \Psi[D_1, D_i]\Big|_{\partial D_1} \simeq
\sqrt {\epsilon_1}.$$ Therefore, we have (iii). \qed

\begin{lem}\label{sum:Hnuh}
We have
$$\left|\int_{\cup_{i=1}^3\p D_i}H\p_\nu h_1\ dS\right|\leq C\sqrt{\epsilon_1}.$$
\end{lem}
\pf Without loss of generality, we assume that
$$\left(-\frac {\epsilon_1} 2, 0\right) \in \p D_1,~\left(\frac {\epsilon_1} 2, 0\right) \in \p D_2,~\mbox{dist}(D_1, D_2) =\epsilon_1~\mbox{and}~(-1,0) \in D_2.$$
We consider $\widetilde{H}$ as follows:
$$\widetilde{H}= H - \p_{2} H(0,0) \frac {x_2}{|\mathbf{x}-(1,0)|^2}.  $$
It follows from the Divergence Theorem that
$$\int_{\p D_1 \cup \p D_2\cup \p D_3} \frac {x_2}{|\mathbf{x}-(1,0)|^2} \p_{\nu} h dS=\int_{\p D_1 \cup \p D_2 \cup \p D_3} \p_{\nu}\left( \frac {x_2}{|\mathbf{x}-(1,0)|^2}\right)  h ds =0,$$
since $\frac {x_2}{|\mathbf{x}-(1,0)|^2} = O (|x|^{-1})$ as
$|x|\rightarrow \infty$. Hence, we have
$$\int_{\cup_{i=1}^3\p D_i}H\p_\nu h_1\ dS = \int_{\p D_1 \cup\p D_2}\widetilde{H}\p_\nu h_1\ dS + \int_{\p D_3}\widetilde{H}\p_\nu h_1\ dS.$$

\par We first consider $\int_{\p D_1 \cup\p D_2}\widetilde{H}\p_\nu h_1\
dS$. By Lemma \ref{lem:A} and \ref{lem:C}, we have $$0>
\partial_{\nu} h_1 \geq C\p_{\nu}\Psi[D_1, D_4]~\mbox{on}~\partial D_1$$
and
$$0 <
\partial_{\nu} h_1 \leq \partial_{\nu}
\Psi[D_1,D_2]~\mbox{on}~\partial D_2.$$ From definition, $\p_{2}
\widetilde{H} =0$. Hence, we can use Lemma 3.2 in \cite{Y2} so that
$$\left|\int_{\p D_1}\widetilde{H}\p_\nu h_1\
dS \right| \leq C \int_{\p D_1}\left|\widetilde{H}\Psi[D_1,
D_4]\right|\ dS \leq C \sqrt {\epsilon_1}$$ and
$$\left|\int_{\p D_2}\widetilde{H}\p_\nu h_1\
dS \right| \leq  \int_{\p D_2}\left|\widetilde{H}\Psi[D_1,
D_2]\right|\ dS \leq C \sqrt {\epsilon_1}.$$ We second consider
$\int_{\p D_3}\widetilde{H}\p_\nu h_1\ dS$. By Lemma \ref{lem:B}, we
can have
$$\left|\int_{\p D_3}\widetilde{H}\p_\nu h_1\
dS\right| \leq C \sqrt {\epsilon_1}.$$ Therefore, we have done it.

\qed

\begin{rem}
We draw attention of readers to the independent work of Bao, Li and
Yin in \cite{Bao_thesis} and \cite{BLY2}.  Bao et al. have shown
that  the  blow-up rate know only for a pair of inclusion is still
valid to the multiple inclusions cases. As a byproduct of our work,
the blow-up rate of the gradient for three inclusions is established
in Theorem \ref{threedomains}.
\end{rem}

\begin{thm}\label{threedomains}
Let $D_1$, $D_2$ and $D_3$ be as assumed in the beginning of Section
\ref{section:threedomains}. Note that $D_2$ is not assumed to be
smaller than the others.
\begin{itemize}
\item[\rm(i)] Optimal upper bounds: For any entire harmonic function $H(x_1,x_2)$, we have the following:
in the narrow region between $D_1\cup D_2$, $$|\nabla u|\leq C\frac{1}{\sqrt{\epsilon_1}},$$
and, in the narrow region between $D_2\cup D_3$, $$|\nabla u|\leq C\frac{1}{\sqrt{\epsilon_2}}.$$
\item[\rm(ii)] Existence of blow-up: Without loss of generality, we assume that
$$\left(-\frac {\epsilon_1} 2, 0\right) \in \p D_1,~\left(\frac {\epsilon_1} 2, 0\right) \in \p D_2 ~\mbox{and}~\mbox{dist}(D_1, D_2) =\epsilon_1.$$
For $H(x_1,x_2)=x_1$, there exist $\mathbf{x}_0$  in the narrow
region between $D_1$ and $D_2$ such that
$$|\nabla u(\mathbf{x}_0)|\geq C\frac{1}{\sqrt{\epsilon_1}},$$
and, similarly, there is a linear function $H (x_1,x_2)$  with
$\mathbf{y}_0$ between $D_2$ and $D_3$ such that
$$|\nabla u(\mathbf{y}_0)|\geq C\frac{1}{\sqrt{\epsilon_2}}.$$
\end{itemize}
\end{thm}

\pf From Subsection \ref{subsec}, we have a representation
\eqref{def:c} for $u$ and the coefficient $c_i$, $i=1,2$, is bounded
independently of $\epsilon_1$ and $\epsilon_2$. Proposition
\ref{prop:nablabound} yields the upper bound of Theorem
\ref{threedomains}.

\par Now, we consider the existence of the blow-up. Using the result
of Subsection \ref{subsec} again, we have a constant $C$ independent
of $\epsilon$ such that
$$\ds\norm{u}_{L^{\infty} (\cup_{i=1}^3 \p D_i)} \leq C \norm{H}_{L^{\infty} (\cup_{i=1}^3 D_i)}.$$
Applying the Green's identity to $\int_{\cup_{i=1}^3\p D_i}u\p_\nu
h_1\ dS$, we have
\begin{align}
\ds&\int_{\cup_{i=1}^3\p D_i}H\p_\nu h_1\ dS = \int_{\cup_{i=1}^3\p
D_i}u\p_\nu h_1\ dS\nonumber\\\nonumber &=-u\bigr|_{\p
D_1}+u\bigr|_{\p D_2}\bigr(\int_{\p D_2}\p_\nu h_1\
dS\bigr)+u\bigr|_{\p D_3}\bigr(\int_{\p D_3}\p_\nu h_1\
dS\bigr)\\\label{eqn:thmC} \ds&=-u\bigr|_{\p D_1}+u\bigr|_{\p
D_2}\bigr(1-\int_{\p D_3}\p_\nu h_1\ dS\bigr)+u\bigr|_{\p
D_3}\bigr(\int_{\p D_3}\p_\nu h_1\ dS\bigr).
\end{align}
By Lemma \ref {lem:B}, we have
\begin{align}
\ds\label{diff:D1D2}u\bigr|_{\p D_2}-u\bigr|_{\p D_1}\leq C
\sqrt{\epsilon_1},\end{align} where the constant $C$ above depends
on $\|H\|_{L^\infty(\cup_{i=1}^3\p D_i)}$. Similarly, we have
\begin{equation}\label{diff:D2D3}u\bigr|_{\p D_3}-u\bigr|_{\p D_2}\leq C \sqrt{\epsilon_2}.\end{equation}
Using \eqref{eqn:thmC} again, we have
\begin{align} u\bigr|_{\p D_2}-u\bigr|_{\p D_1} + O (\sqrt {\epsilon_1
\epsilon_2})&= \int_{\cup_{i=1}^3\p D_i}H\p_\nu h_1\ dS
\notag\\&\geq\int_{\p D_1}H\p_\nu h_1\ dS. \label{eq:4}
\end{align}
The last inequality can be derived from the fact that $H >0$ on $\p
D_2 \cup \p D_3$.

 To get the last inequality above, we
took an advantage of $H=x_1$. By \eqref{eq:3}, $$  \partial_{\nu} h
\leq \left(\frac{ h_1\Bigr|_{\p D_1}-h_1\Bigr|_{\p D_2\cup\p
D_3}}{\Psi[D_1,D_2]\Bigr|_{\p D_1}-\Psi[D_1,D_2]\Bigr|_{\p
D_2}}\right) \partial_{\nu} \Psi[D_1,D_2]<0~\mbox{on}~\partial D_1.
$$ By (iii) in Proposition \ref{prop:nablabound}, we have
$$ h_1\Bigr|_{\p D_1}-h_1\Bigr|_{\p D_2\cup\p
D_3} \backsimeq \sqrt {\epsilon_1}$$and$${\Psi[D_1,D_2]\Bigr|_{\p
D_1}-\Psi[D_1,D_2]\Bigr|_{\p D_2}}\backsimeq \sqrt {\epsilon_1}.$$
The inequality \eqref{eq:4} implies
$$ u\Bigr|_{\p D_2}-u\Bigr|_{\p D_1} \gtrsim \sqrt {\epsilon}.$$
By the Mean Value Theorem, we have the desirable lower bound in the
narrow region between $D_1$ and $D_2$. Similarly, we can also obtain
the other lower bound.
 \qed

\subsection{Proof of Theorem \ref{thm:B}}
We derive the optimal bounds of the gradient of the solution to \eqref{eq:conductivity}, when there are adjacent three disks:
 \begin{equation}\label{def:Bi}D_l = B_{r_l}(\mathbf{c}_l),\
l=1,2,3,\end{equation}  where $\mathbf{c}_1=(-r_1 - \frac
{\epsilon_1} 2, 0)$, $\mathbf{c}_2=(r_2+ \frac {\epsilon_1} 2,0)$
and $\mathbf{c}_3=(r_3+r_2+\frac {\epsilon_1} 2 + \epsilon_2,0)$. As
defined before, $h_1 = \Psi [D_1, (D_2\cup D_3)]$. Let $w_1 =\Psi
[D_1, D_2] $.

\par We begin the proof by showing  that
\be w_1 \big|_{\partial D_2} -  w_1 \big|_{\partial D_1} \simeq h_1
\big|_{\partial D_2} -  h_1 \big|_{\partial D_1} .\label{eq:5}\ee By
the monotonic property of Lemma \ref{lem:mono}, we have
$$ h_1 \big|_{\partial D_2} -  h_1 \big|_{\partial D_1} \leq
w_1 \big|_{\partial D_2} -  w_1 \big|_{\partial D_1}.$$ Considering
$$h_1 - \left({\frac {h_1 \big|_{\partial D_2} -  h_1 \big|_{\partial D_1} }{w_1 \big|_{\partial D_2} -  w_1 \big|_{\partial D_1}}}\right)w_1,$$
we can obtain, from the Hopf's Lemma,
$$\int_{\p D_2} \p_{\nu}h_1 dS \leq \left({\frac {h_1 \big|_{\partial D_2} -  h_1 \big|_{\partial D_1} }{w_1 \big|_{\partial D_2} -  w_1 \big|_{\partial D_1}}}\right)\int_{\p D_2} \p_{\nu}w_1 dS.$$
By Lemma \ref{lem:B}, we have
$$ \int_{\p D_3} \p_{\nu}h_1 dS = O(\sqrt {\epsilon_1}).$$
Since $\int_{\p D_2 \cup \p D_3} \p_{\nu}h_1 dS =1$, we have
$$ \left(w_1 \big|_{\partial D_2} -  w_1 \big|_{\partial D_1}\right)
(1+ O(\sqrt {\epsilon_1})) \leq h_1 \big|_{\partial D_2} -  h_1
\big|_{\partial D_1} .$$ Therefore, we can obtain \eqref{eq:5}.
Owing to the estimate for $w_1 \big|_{\partial D_2} -  w_1
\big|_{\partial D_1}$ in Lemma \ref{lemma:twodisks}, we have

\be h_1 \big|_{\partial D_2} -  h_1 \big|_{\partial D_1} \simeq
\sqrt {\frac {r_1 + r_2}{r_1 r_2}} \sqrt
{\epsilon_1}.\label{eq:6}\ee

\par Let $w_2 =\Psi[D_1, D_3] $. Considering

$$h_1 - \left({\frac {h_1 \big|_{\partial D_3} -  h_1 \big|_{\partial D_1} }{w_2 \big|_{\partial D_3} -  w_2 \big|_{\partial D_1}}}\right)w_2,$$
 from the Hopf's Lemma, we obtain
$$\p_{\nu} h_1 \leq  \left({\frac {h_1 \big|_{\partial D_3} -  h_1 \big|_{\partial D_1} }{w_2 \big|_{\partial D_3} -  w_2 \big|_{\partial D_1}}}\right)\p_{\nu}w_2 \leq 0~\mbox{on}~\partial D_1.$$
Here, we estimate the coefficient in the right hand side. Note that
$h_1 \big|_{\partial D_2} =  h_1 \big|_{\partial D_3} $. Thus, we
have $$h_1 \big|_{\partial D_3} -  h_1 \big|_{\partial D_1} \simeq
\sqrt {\frac {r_1 + r_2}{r_1 r_2}} \sqrt {\epsilon_1}.$$ Since $r_2
\ll r_1$ and $r_2 \ll r_3$, we also have
$$w_2 \big|_{\partial D_3} -  w_2 \big|_{\partial D_1} \simeq
\sqrt {\frac {r_1 + r_3}{r_1 r_3}} \sqrt {r_2}.$$ This implies that
$$\p_{\nu} h_1 \lesssim  \frac {\sqrt {\frac {r_1 + r_2}{r_1 r_2}} \sqrt {\epsilon_1}}  {\sqrt {\frac {r_1 + r_3}{r_1 r_3}} \sqrt {r_2}} \p_{\nu} w_2 \leq 0~\mbox{on}~\p D_1.$$
Therefore, we have
\begin{align}
\int_{\p D_1} H \p_{\nu} h_1 dS &\gtrsim \sqrt {\frac {r_1 +
r_2}{r_1 + r_3}} \frac {\sqrt {r_3}}{r_2} \sqrt {\epsilon_1}
\int_{\partial D_1} H \p_{\nu} w_2 dS\notag\\&\gtrsim
 \sqrt {\frac {r_1 +
r_2}{r_1 + r_3}} \frac {\sqrt {r_3}}{r_2} \sqrt {\epsilon_1}
\sqrt{\frac {r_1 r_3}{r_1 + r_3}} \sqrt {r_2}\notag\\
&\geq \frac {r_1 r_3}{r_1 + r_3} \frac 1 {\sqrt {r_2}}
\sqrt{\epsilon_1} \geq 0 \label{eq:8}.
\end{align}
Owing to \eqref{eqn:thmC}, \eqref{diff:D1D2} and \eqref{diff:D2D3},
we have
\begin{align}
\int_{\cup_{i=1}^3\p D_i}u \p_\nu h_1\ dS =&\left(1- \int_{\p
D_3}\p_\nu h_1\ dS\right)\left( u\bigr|_{\p D_2} -u\bigr|_{\p
D_1}\right)\notag\\&~~~~~+  \left( \int_{\p D_3}\p_\nu h_1\ dS
\right) \left( u\bigr|_{\p D_3} -u\bigr|_{\p
D_1}\right) \notag\\
=& \left(1- O (\sqrt{\epsilon_1})\right)\left( u\bigr|_{\p D_2}
-u\bigr|_{\p D_1}\right)+  O (\sqrt{\epsilon_1})\left( O
(\sqrt{\epsilon_1}) + O (\sqrt{\epsilon_2})\right). \notag
\end{align}
Therefore, we have
\begin{align}
u\bigr|_{\p D_2} -u\bigr|_{\p D_1} & \geq  \frac 1 2
\int_{\cup_{i=1}^3\p D_i}u \p_\nu h_1\ dS + O
(\sqrt{\epsilon_1})\left( O (\sqrt{\epsilon_1}) + O
(\sqrt{\epsilon_2})\right) \notag\\
&=  \frac 1 2 \int_{\cup_{i=1}^3\p D_i} H \p_\nu h_1\ dS + O
(\sqrt{\epsilon_1})\left( O (\sqrt{\epsilon_1}) + O
(\sqrt{\epsilon_2})\right) \notag\\
&\geq \frac 1 2 \int_{\p D_1} H \p_\nu h_1\ dS + O
(\sqrt{\epsilon_1})\left( O (\sqrt{\epsilon_1}) + O
(\sqrt{\epsilon_2})\right) \notag\\
&\geq C {\frac {r_1 r_3}{r_1 + r_3}} {\frac 1 {\sqrt {r_2}}} \sqrt
{\epsilon_1} + O (\sqrt{\epsilon_1})\left( O (\sqrt{\epsilon_1}) + O
(\sqrt{\epsilon_2})\right).\notag
\end{align}
Therefore, we have completed the proof. \qed

\subsection{Proof of Theorem \ref{thm:D}}
We pursuit the proof of Theorem \ref{thm:B}, taking an advantage of
the monotonic property of Lemma \ref{lem:mono}. The domains $D_1$,
$D_2$ and $D_3$ are as assumed in Theorem \ref{thm:D}. As assumed
before, $h_1 =\Psi [D_1, (D_2\cup D_3)]  $. Let $w_1 =\Psi [D_1,
D_2] $. By the same way as Theorem \ref{thm:B}, we have $$ w_1
\big|_{\partial D_2} -  w_1 \big|_{\partial D_1} \simeq h_1
\big|_{\partial D_2} - h_1 \big|_{\partial D_1} .$$ Here, we use the
monotonic property of Lemma \ref{lem:mono} to estimate the
difference between domains. Choosing two pairs of proper disks
containing $D_1$ and $D_2$, and contained $D_1$ and $D_2$,
respectively, we can obtain $$ h_1 \big|_{\partial D_2} - h_1
\big|_{\partial D_1} \simeq \sqrt{\frac {\epsilon_1}{ { r_2}}}$$
under the assumption that $r_2$ is small.
\par  Let $w_2 =\Psi[D_1, D_3] $. Choosing two
pairs of proper disks containing $D_1$ and $D_3$, and contained
$D_1$ and $D_3$, respectively, Then, we have
$$w_2 \big|_{\partial D_3} -  w_2 \big|_{\partial D_1} \simeq
\sqrt {r_2}.$$ By the same argument as Theorem \ref{thm:B}, we have
$$
\int_{\p D_1} H \p_{\nu} h_1 dS \gtrsim \sqrt {\frac{\epsilon_1}
{r_1}} \geq 0.
$$
Note that $D_1 \subset \mathbb{R}_{-}\times \mathbb {R}$ and
$D_2\cup D_3 \subset \mathbb{R}_{+}\times \mathbb {R}.$ Continuing
to follow the proof of Theorem \ref{thm:B}, we can obtain
$$
u\bigr|_{\p D_2} -u\bigr|_{\p D_1}  \geq C  \sqrt {\frac
{\epsilon_1}{r_1}} + O (\sqrt{\epsilon_1})\left( O
(\sqrt{\epsilon_1}) + O (\sqrt{\epsilon_2})\right).\notag
$$
Therefore, we have done the proof.

\subsection {Derivation for the optimal upper bounds}\label{upp}
We consider the optimal upper bounds presented in Theorem
\ref{thm:A}, \ref{thm:B}, \ref{thm:C} and \ref{thm:D}.  These proofs
have essential thing in common. In this respect, we prove only the
optimal upper bound presented in Theorem \ref{thm:B}. As have
assumed them before, we set \begin{align*} h_1 &=\Psi [D_1, D2 \cup D_3]\\
 h_2 &=\Psi [D_1, D_2]\\
  h_3 &=\Psi [D_1, D_3]\\
   h_4 &=\Psi [D_1, D_4].
\end{align*}
Here, the domain $D_4$ is given in Lemma \ref{lem:C}, which is a
disk containing $D_2$ and $D_3$ with
$$\mbox{dist}(D_1,D_4) = \mbox{dist}(D_1,D_2),$$
and the diameter of $D_4$ is in proportion as $r_3$, because $r_2$
is sufficiently small. Then, we compare $ h_1$ with $h_2$, $h_3$ and
$h_4$. The proof of Lemma \ref{lem:A} contains
$$0 \leq \partial_{\nu} h_1 \leq \left(\frac{ h_1\Bigr|_{\p (D_2 \cup D_3)}-h_1\Bigr|_{\p D_1}}{h_2\Bigr|_{\p D_2 }-h_2\Bigr|_{\p
D_1}}\right) \partial_{\nu}h_2~\mbox{on}~\p D_2,
$$ the proof of Lemma \ref{lem:B} yields
$$0 \leq \partial_{\nu} h_1 \leq \left(\frac{ h_1\Bigr|_{\p (D_2 \cup D_3)}-h_1\Bigr|_{\p D_1}}{h_3\Bigr|_{\p D_3 }-h_3\Bigr|_{\p
D_1}}\right) \partial_{\nu}h_3~\mbox{on}~\p D_3
$$ and the proof of Lemma \ref{lem:C} implies
$$0 \leq -\partial_{\nu} h_1 \leq -\left(\frac{ h_1\Bigr|_{\p (D_2 \cup D_3)}-h_1\Bigr|_{\p D_1}}{h_4\Bigr|_{\p D_4 }-h_4\Bigr|_{\p
D_1}}\right) \partial_{\nu}h_3~\mbox{on}~\p D_3.
$$ In the same way as Lemma \ref {sum:Hnuh}, we can consider $\widetilde{H}$ by choosing the point in $D_3$. In this respect, without any loss of generality,
we can assume that
$$\partial_{x_2} H (0,0) = 0.$$ The reason why we assumed above is
because the integration representation for the potential difference
is not good enough, refer to \cite{Y2}. The geometrical assumption
of Case (B) implies that $D_1$ and $D_2 \cup D_3$ are separated by
$x_1 = 0$ and they are approaching to $(0,0)$.
\par Therefore, by the proof of Theorem \ref{thm:B} and Lemma \ref{lem++}, we have
\begin{align*}
\left| u\bigr|_{\p D_2} -u\bigr|_{\p D_1}\right| &+ O
(\sqrt{\epsilon_1})\left( O
(\sqrt{\epsilon_1}) + O (\sqrt{\epsilon_2})\right)\\
&= \left| \int_{\partial (\cup_{i=1} ^3 D_i)} H \partial_{\nu} h_1
dS \right|\\
&\lesssim \left(\frac{ h_1\Bigr|_{\p (D_2 \cup D_3)}-h_1\Bigr|_{\p
D_1}}{h_2\Bigr|_{\p D_2 }-h_2\Bigr|_{\p D_1}}\right)
\sqrt{\frac{r_1r_2}{r_1+r_2}\epsilon_1}\\
&~~+ \left(\frac{ h_1\Bigr|_{\p (D_2 \cup D_3)}-h_1\Bigr|_{\p
D_1}}{h_3\Bigr|_{\p D_3 }-h_3\Bigr|_{\p D_1}}\right)
\sqrt{\frac{r_1r_3}{r_1+r_3}r_2}\\&~~ +\left(\frac{ h_1\Bigr|_{\p
(D_2 \cup D_3)}-h_1\Bigr|_{\p D_1}}{h_4\Bigr|_{\p D_4
}-h_4\Bigr|_{\p D_1}}\right) \sqrt{\frac{r_1r_3}{r_1+r_3}\epsilon_1}
\end{align*} and
$$ { h_1\Bigr|_{\p (D_2 \cup D_3)}-h_1\Bigr|_{\p
D_1}} \thickapprox{h_2\Bigr|_{\p D_2 }-h_2\Bigr|_{\p D_1}}.$$ Here,
note that the radius of $D_4$ can be choosen between $\frac 3 2 r_3$
and $2 r_3$. Lemma \ref{lemma:twodisks} implies that $$\left|
u\bigr|_{\p D_2} -u\bigr|_{\p D_1}\right| \lesssim \frac {r_1 r_3}{
r_1 + r_3 } \frac 1 {\sqrt {r_2}} {\sqrt {\epsilon_1}}.$$ Therefore,
we establish the optimal upper bound for $\left| u\bigr|_{\p D_2}
-u\bigr|_{\p D_1}\right| $.

\par Based on this, the optimal upper bound on the
gradient of $u$ in the narrow region be obtained. Here, the main
idea to get the gradient estimate from the potential difference has
already been presented by Bao et al. (Theorem 1.3, Lemma 2.2 and 2.3
in \cite{BLY}), and has been modified to fit our problem by Lim and
Yun in \cite{LY}. Thurs, we give a brief description on the method.
We choose a large domain $D_0$ containing $D_1$, $D_2$ and $D_3$,
where $\p D_0$ is at a sufficient distance from $D_1$, $D_2$ and
$D_3$. Then, $u$ can be decomposed as follows:
$$u=C_0 + v_0 + C_1 v_1  + C_3 v_3$$
where  for $i=0,1,3$, $v_i$ is a harmonic function in $D_0 \setminus
(D_1\cup D_2 \cup D_3)$ with the boundary data
$$v_i = \delta_{0j} ~\mbox{on}~\partial D_j~\mbox{for}~i=1,~3$$
and
$$v_0 = \delta_{ij} u  ~\mbox{on}~\partial D_j$$
for any $j=0,1,2,3$. Thus, the constants $C_1$ and $C_3$ keep
$$ |C_1|  \lesssim \frac {r_1 r_3}{
r_1 + r_3 } \frac 1 {\sqrt {r_2}}{\sqrt {\epsilon_1}}$$ and
$$ |C_3|  \lesssim \frac {r_1 r_3}{
r_1 + r_3 } \frac 1 {\sqrt {r_2}} {\sqrt {\epsilon_2}}.$$

\par To estimate $\nabla v_0$, we consider a harmonic function $\rho$
in $D_0 \setminus (D_1\cup D_2 \cup D_3)$ with the boundary data
$$\rho = \delta_{0j}~ \mbox{on}~\partial D_j$$
for  any $j=0,1,2,3$. By comparing with the harmonic function
$\rho_{i}$ in $D_0 \setminus D_i$ with $\rho_{i} = 0$ on $\p D_i$
and $\rho_{i} = 1$ on $\p D_0$, the Hopf's Lemma yields

\begin{align*}&\norm { \nabla \rho }_{L^{\infty} (D_0 \setminus (D_1\cup D_2
\cup D_3 ))} \\ &\leq \max \{\norm { \nabla \rho_1 }_{L^{\infty} (\p
D_1 )}, \norm { \nabla \rho_2 }_{L^{\infty} (\p D_2 )}, \norm {
\nabla \rho_3 }_{L^{\infty} (\p D_3 )},\norm { \nabla \rho
}_{L^{\infty} (\p D_0 )}\}< C.\end{align*} Applying the Hopf's Lemma
again, we can have that the gradient of $v_0$ is bounded independent
of $\epsilon_1$, refer to Lemma 2.2 in \cite{BLY}.

\par We estimate $C_1 \nabla v_1$ in the narrow region between $D_1$ and
$D_2$. Since $v_1$ is constat on the boundaries and the boundaries
is smooth enough in the narrow region, the proof of Lemma 4.3
implies that $v_1$ can be extend into the interior areas of $D_1$
and $D_2$ by the distance almost $\epsilon$ from the boundaries in
the narrow region, independently of $r_1$ and $r_2$. By the gradient
estimate for harmonic functions allows
$$ |C_1 \nabla v_1 | \lesssim \frac {r_1 r_3}{
r_1 + r_3 } \frac 1 {\sqrt {r_2}}\frac 1 {\sqrt {\epsilon_1}}$$ in
the narrow region between $D_1$ and $D_2$. Note that the inequality
above is a local property independent of choosing $D_0$.

\par Now, we consider $ C_3 \nabla v_3$ in the narrow region between $D_1$ and
$D_2$. Let $\widetilde{\rho}$ be a harmonic function in $D_0
\setminus (D_2 \cup D_3)$ with the boundary data $$ \widetilde{\rho}
= 0 ~\mbox{on}~\partial (D_0 \cup D_2) ~\mbox{and}~\widetilde{\rho}
= 1 ~\mbox{on}~\partial D_3.
$$ By the maximum principle, we have
$$0\leq   v_3 \leq \widetilde{\rho}~\mbox{in}~D_0
\setminus (D_1\cup D_2 \cup D_3)$$ Considering the standard estimate
for $\Psi[D_2,D_3]$, we can obtain
$$ |C_3 v_3| \leq C \sqrt {\epsilon_2}.$$
Similarly to the estimate for $C_1 \nabla v_1$, the gradient
estimate for harmonic functions yields
$$ |C_3 \nabla v_3 | \lesssim \sqrt  {\frac {\epsilon_2} {{\epsilon_1}}}$$
in the narrow region between $D_1$ and $D_2$. \par Therefore, we can
obtain the desirable upper bound. Here, it is noteworthy that the
upper bound is dominated only by the estimate for $C_2 \nabla v_1$,
which is independent of choosing $D_0$. In this respect, the
constant C of the upper bound in Theorem \ref{thm:B} is independent
of $r_1$, $r_2$, $r_3$, $\epsilon_1$ and $\epsilon_2$.
 \qed

\section*{Acknowledgements}
The authors would like to express his gratitude to \textrm{Professor
Hyeonbae Kang}, who suggested the original problem studied in this
paper.  The authors is also grateful to \textrm{Professor YanYan Li}
for his concern with their subject and suggestions. The second named
author would also like to express his thanks to \textrm{Professor
Gang Bao}, and gratefully acknowledges his hospitality during the
visiting period at Michigan State University.


\begin{thebibliography}{}\label{sec:TeXbooks}


\bibitem{A}{\sc L. Ahlfors}, {Complex Analysis}, Third ed., McGraw-Hill, New
York, 1979.



\bibitem{ADKL} {\sc Ammari H, Dassios G, Kang H, and Lim M}, {\em Estimates for the electric field in the presence of adjacent perfectly conducting spheres}, {Quat. Appl. Math.}, { 65} (2007), pp.~339--355

\bibitem{AKL} {\sc H. Ammari, H. Kang, and M. Lim}, {\em Gradient estimates for solutions to the conductivity problem}, Math. Ann., 332(2) (2005), pp.~277--286.

\bibitem{AKLLL} {\sc H. Ammari, H. Kang, H. Lee, J. Lee, and M. Lim}, {\em Optimal bounds on the gradient of solutions to conductivity problems}, J. Math. Pures Appl., 88 (2007), pp.~307--324.

\bibitem{Bao_thesis}{\sc S. Bao}, {\em Gradient estimates for the conductivity problems}, thesis, Rutgers University.


\bibitem{BC}{\sc  B. Budiansky and G. F. Carrier}, {\em High shear stresses in stiff fiber composites}, J. Appl. Mech., 51 (1984), pp. 733-735.

\bibitem{BLY} {\sc E. Bao, Y.Y. Li and B. Yin} {\em Gredient estimats for the conductivity problem}, Arch. Rational Mech. Anal, to appear.

\bibitem{BLY2} {E. Bao, Y.Y. Li and B. Yin} {\em Gradient estimates for the perfect and
insulated conductivity problem and elliptic systems} (in
preparation)

\bibitem{BV} {\sc E. Bonnetier and M. Vogelius}, {\em An Elliptic regularity result for a composite medium with "touching" fibers of circular
cross-section}, SIAM J. Math. Anal., 31, No 3 (2000), pp.~651--677.




%\bibitem{J} {\sc J.D. Jackson}, {\em Classical
%Electrodynamics}, Third ed., Wiley, New York, NY, 1999.
%


\bibitem{K} {\sc J.B. Keller}, {\em Stresses in narrow regions}, Trans. ASME
J. Appl. Mech., 60 (1993), pp. ~1054--1056.



\bibitem{LN} {\sc Y.Y. Li and M. Nirenberg}, {\em Estimates for ellliptic system from composite material}, Comm. Pure Appl. Math., LVI (2003), pp.~892--925.


\bibitem{LV} {\sc Y.Y. Li and M. Vogelius}, {\em Gradient estimates for solution to divergence form elliptic equation with discontinuous coefficients}, Arch. Rational Mech. Anal., 153 (2000),
pp.~91--151.

\bibitem{LY} {\sc M. Lim and K. Yun}, {\em Blow-up of Electric Fields between Closely Spaced Spherical Perfect Conductors}, submitted.

\bibitem{Y} {\sc K. Yun},  {\em  Estimates for electric fields blown up between closely adjacent conductors with arbitrary shape}, { SIAM J. Appl. Math.}, 67, No 3 (2007), pp.~714--730.
\bibitem{Y2} {\sc K. Yun}, {\em Optimal bound on high stresses occurring between stiff fibers with arbitrary shaped cross sections}, J. Math. Anal. Appl. 350, (2009), pp. 306-312
\end{thebibliography}
\end{document}